\documentclass[a4paper,10pt,twoside]{amsart}

\usepackage[
  marginratio=1:1,
  height=597pt,
  width=460pt,
  tmargin=117pt
]{geometry}

\usepackage{amsmath,amssymb,amsthm,mathtools}
\usepackage{amsfonts}
\usepackage{mathrsfs}
\usepackage{tikz}
\usepackage{tikz-cd}
\usepackage{hyperref}
\usepackage{enumitem}
\usepackage{graphicx}

\newcommand{\M}{\widetilde{D}_0^\#(\infty)}

\DeclareMathOperator{\Ind}{Ind}
\DeclareMathOperator{\End}{End}
\DeclareMathOperator{\soc}{soc}

\DeclareMathOperator{\Ext}{Ext}
\DeclareMathOperator{\Hom}{Hom}
\DeclareMathOperator{\Sym}{Sym}
\DeclareMathOperator{\GL}{GL}
\newcommand{\Det}{\mathrm{det}}

\usepackage{etoolbox}

\expandafter\patchcmd\csname\string\proof\endcsname
  {\normalparindent}{0pt }{}{}

\makeatletter
\patchcmd{\@thm}
  {\thm@headfont{\scshape}}
  {\thm@headfont{\scshape\bfseries}}
  {}{}

\patchcmd{\@thm}
  {\thm@notefont{\fontseries\mddefault\upshape}}
  {}
  {}{}
\makeatother

\newtheorem{theorem}{Theorem}[section]
\newtheorem{proposition}[theorem]{Proposition}
\newtheorem{lemma}[theorem]{Lemma}

\newtheorem{remark}[theorem]{Remark}

\theoremstyle{definition}

\title[Nonadmissible irreducible representations]{On the existence of nonadmissible irreducible smooth representations of $p$-adic \(\mathrm{GL}_n\)}

\author{DANIEL LE, MANSIMAR SINGH}

\date{}

\numberwithin{equation}{section}

\begin{document}

\begin{abstract}
Let \(p>7\) be a prime, let \(n>1\), and let \(F\) be a nonarchimedean local field of residue characteristic \(p\) not equal to $\mathbb Q_p$. For every extension \(E\) of the residue field of \(F\), we construct an irreducible smooth \(E\)-representation of \(\mathrm{GL}_n(F)\) which is not admissible. When $E$ is countable, we also construct such representations without a central character and with nonscalar endomorphisms. In the case where the residue field is \(\mathbb F_p\), the new ingredient is a construction of Breuil--Pa\v{s}k\={u}nas diagrams which are not invariant under the first principal congruence subgroup.
\end{abstract}

\maketitle

\section{Introduction}

Let \(F\) be a nonarchimedean local field of residue characteristic \(p\). A basic problem in the smooth representation theory of \(\mathrm{GL}_n(F)\) is to determine which finiteness properties follow from irreducibility. When the coefficient field has a characteristic different from \(p\), irreducible smooth representations are admissible \cite{Vig96}. In defining characteristic, the situation is substantially less rigid.

\begin{theorem}\label{thm:main}
Let \(n>1\), and let \(p>7\) be a prime. Let \(F\neq\mathbb Q_p\) be a nonarchimedean local field of residue characteristic \(p\).

\begin{enumerate}
\item Let \(E\) be an extension of the residue field of \(F\) of countable cardinality. Then there exists an irreducible smooth \(E\)-representation of \(\mathrm{GL}_n(F)\) which does not admit a central character.

\item Let \(E\) be an extension of the residue field of \(F\), and let
\(\widetilde E\) be a division \(E\)-algebra of countable \(E\)-dimension
with center \(\widetilde C\). Then there exists an irreducible
nonadmissible smooth \(E\)-representation of \(\mathrm{GL}_n(F)\) whose
endomorphism algebra is isomorphic to \(\widetilde E\) as an \(E\)-algebra.
Moreover, it may be chosen
\begin{itemize}
\item with a Hecke eigenvalue and a central character;
\item without a central character and without a Hecke eigenvalue, if
\(E\subsetneq\widetilde C\);
\item with a central character but without a Hecke eigenvalue, if
\(n>3\) and \(E\subsetneq\widetilde C\).
\end{itemize}
\end{enumerate}
\end{theorem}

\begin{remark}
    When $n=2$, Theorem~\ref{thm:main} also holds for $p=7$.
\end{remark}

The case in which the residue field of \(F\) is different from
\(\mathbb F_p\) is contained in \cite{Le24}, building on
\cite{Le19,GS20, GLS23}. It remains to prove Theorem~\ref{thm:main} when the
residue field of \(F\) is \(\mathbb F_p\) and \(F\neq\mathbb Q_p\). This is
the case treated in this paper.

We briefly explain why this remaining case requires a new construction. The
basic tool in the previous work is a diagram in the sense introduced by
Breuil and Pa\v{s}k\={u}nas \cite{BP12}. The diagrams used in
\cite{Le19,GLS23,Le24} are invariant under the first principal congruence
subgroup of \(\mathrm{GL}_2(\mathcal O_F)\). At this level, the category of
such diagrams depends only on the residue field: if two local fields have the
same residue field, then their corresponding categories of diagrams are
equivalent. Thus, a construction using only such diagrams over a field with
residue field \(\mathbb F_p\) would yield the same kind of construction over
\(\mathbb Q_p\).

This is a genuine obstruction to the previous method. Indeed, in rank two,
every irreducible smooth \(\overline{\mathbb F}_p\)-representation of
\(\mathrm{GL}_2(\mathbb Q_p)\) is admissible and has a central character
\cite{BL94,Bre03,Ber12}. Thus, in the remaining case, a construction of
irreducible representations without a central character must use information
that is invisible at the first congruence level.

Our construction passes to the second congruence level. The resulting
diagrams are not invariant under the first principal congruence subgroup, and
therefore retain information that is lost modulo the first congruence
subgroup. This is reflected in the diamond-shaped Hasse diagram of the cosocle filtration of
Proposition~\ref{diamond-prop}. When \(F=\mathbb Q_p\), the corresponding
diamond is replaced by the Hasse diagram of
Remark~\ref{rmk:Qp}. This diamond-shaped filtration is the key new feature of
the construction and provides the flexibility needed to prove
Theorem~\ref{thm:main}.

\subsection{Related work}

In independent recent work, Huang \cite{Huang} constructs, for \(p\geq 5\), nonadmissible smooth irreducible \(\overline{\mathbb F}_p\)-representations of \(\mathrm{GL}_2(F)\) when \(F/\mathbb Q_p\) is a proper totally ramified finite extension. While, unsurprisingly, both constructions use diagrams that are not \(K_1\)-invariant, the constructions are otherwise fairly distinct. 
We further obtain a result for the equal characteristic case, some results beyond nonadmissibility, and results for \(\mathrm{GL}_n(F)\) for all \(n>1\) using parabolic induction, though we expect similar results could be obtained using Huang's methods. 

\subsection{AI Statement}

ChatGPT was used for language editing, improving the exposition, and assisting in checking mathematical arguments for possible errors or points requiring clarification. The mathematical content was developed and verified by the authors.

\subsection{Acknowledgments}

The results in Section \ref{sec:diamond} were developed in a collaboration between Bao Le Hung, Zhixiang Wu, and D.L. 
We thank them for allowing us to include this here. 
We thank Dayu Huang for sharing his work \cite{Huang} with us and for his interest in our work. 
We thank Yongquan Hu for comments about Remark \ref{rmk:Qp}. 
D.L. was supported by the National Science Foundation under agreement DMS-2302623. 

\subsection{Notation and conventions}

Let \(p>5\) be a prime, and let \(F\) be a nonarchimedean local field whose residue field is \(\mathbb F_p\) such that \(F\neq\mathbb Q_p\). Let
\(\mathcal O_F\) denote its ring of integers, fix a uniformizer
\(\varpi\in\mathcal O_F\), and let \(E\) be an arbitrary field extension of
\(\mathbb F_p\). The only feature of the hypothesis on \(F\) used below is that
\(\mathcal O_F/\varpi^2\) is the ring of dual numbers over \(\mathbb F_p\).
We fix the identification $\mathcal O_F/\varpi^2\cong \mathbb F_p[\varepsilon]/(\varepsilon^2),
$ under which the image of \(\varpi\) is \(\varepsilon\).

Set \(G:=\mathrm{GL}_2(F)\) and \(K:=\mathrm{GL}_2(\mathcal O_F)\), and let
\(Z\subseteq G\) denote the center. Let \(I\subseteq K\) be the upper triangular Iwahori subgroup, \(I_1\subseteq I\) the normal pro-\(p\) Sylow
subgroup, and let $\Pi:=\bigl(\begin{smallmatrix} 0 & 1 \\ \varpi & 0 \end{smallmatrix}\bigr) \in N_G(I)$. For an \(I\)-module \(V\), let \(V^{(\Pi)}\) denote the \(\Pi\)-twisted \(I\)-module \(\{v^{(\Pi)}\mid v\in V\}\), with \(I\)-action given by \(g\cdot v^{(\Pi)}=(\Pi g\Pi^{-1}\cdot v)^{(\Pi)}\) for \(g\in I\) and \( v\in V\).

Let \(K_1:=1+\varpi M_2(\mathcal O_F)\) and
\(K_2:=1+\varpi^2M_2(\mathcal O_F)\). Reduction modulo \(\varpi^2\)
identifies
\[
K/K_2\cong
\mathrm{GL}_2(\mathcal O_F/\varpi^2)\cong
\mathrm{GL}_2(\mathbb F_p[\varepsilon]/(\varepsilon^2)).
\]
Under this identification, \(K_1/K_2\) corresponds to
\(1+\varepsilon M_2(\mathbb F_p)\), and
\(\Gamma:=K/K_1\cong\mathrm{GL}_2(\mathbb F_p)\) is identified with the image
of the inclusion \(\mathrm{GL}_2(\mathbb F_p)\hookrightarrow
\mathrm{GL}_2(\mathbb F_p[\varepsilon]/(\varepsilon^2))
\) as matrices with zero \(\varepsilon\)-coefficient. This gives a section
\begin{equation}\label{section}
\Gamma\hookrightarrow K/K_2
\end{equation}
of the projection \(K/K_2\twoheadrightarrow\Gamma\), and hence a semidirect
product decomposition \(K/K_2\cong (K_1/K_2)\rtimes\Gamma\).

Any simple smooth \(E[K]\)-module is obtained by inflation from \( \Gamma\), is defined over \(\mathbb{F}_p\), and is absolutely irreducible; see, for example, \cite[\S~1]{BL94}. Moreover, taking \(\Gamma\)-invariants commutes with extension of scalars from \(\mathbb{F}_p\) to \(E\), as does \(\Hom_{\Gamma}(V,-)\) for finite-dimensional $V$. In particular, formation of socle and cosocle for finite-dimensional representations commutes with base change.

Let \(\rho\) be a Serre weight, i.e., an irreducible smooth
\(K\)-representation over \(E\). By the standard classification of weights
\cite[Proposition~2.2]{Her11a}, Serre weights are in bijection with
\[
\{(a,b)\in\mathbb Z^2\mid 0\leq a-b\leq p-1\}
/\Delta((p-1)\mathbb Z),
\]
where \(\Delta\) denotes the diagonal embedding of \((p-1)\mathbb Z\) into
\(\mathbb Z^2\). The class of \((a,b)\) corresponds to
\[
F(a,b):=\mathrm{Sym}^{a-b}E^2\otimes\Det^b .
\]
A Serre weight \(F(a,b)\) is called \textit{generic} if
\(a-b\notin \{ 0,p-1\}\). Since \(\rho^{I_1}\) is
one-dimensional, we write \(\chi(\rho)\) for the character of \(I\) acting on
this line. The map \(\rho\mapsto\chi(\rho)\) gives a bijection from the set
of generic Serre weights to the set of \(I\)-characters \(\chi\) such that
\(\chi\neq\chi^{(\Pi)}\). For a generic Serre weight \(\rho\), we denote by
\(\rho^{[s]}\) the unique generic Serre weight with
\(\chi(\rho^{[s]})=\chi(\rho)^{(\Pi)}\).

When \(\rho=F(a,b)\) with \(a-b\geq 1\), let \(V_2(\rho)\) denote the
two-dimensional step of the unique \(I\)-stable full flag in \(\rho|_I\). The
representation \(V_2(\rho)\) fits into the non-split exact sequence
\begin{equation}\label{V2}
0\longrightarrow \chi(\rho)\longrightarrow V_2(\rho)
\longrightarrow \chi(\rho)\alpha^{-1}\longrightarrow 0,
\end{equation}
where \(\alpha:I\to E^\times\) is the character
\[
\alpha\!\left(\begin{pmatrix} x&y\\ z&w\end{pmatrix}\right)
=\overline{xw^{-1}},
\]
with \(\overline{xw^{-1}}\) denoting reduction modulo \(\varpi\).

\section{A diamond module}\label{sec:diamond}

Fix a Serre weight \(\sigma:=F(a,b)\) with \(r:=a-b\) satisfying
\(3\leq r\leq p-2\), so that \(\sigma^{[s]}=F(p-1+b,a)\). Also fix
\(\chi:=\chi(\sigma)\), and let \(\tau\) be the unique Serre weight with
\(\chi(\tau)=\chi\alpha^{-1}\). Explicitly,
\[
\tau=F(a-1,b+1)
\quad \text{and}\quad
\tau^{[s]}=F(p+b,a-1).
\]

Proposition~\ref{diamond-prop} is the first step in our construction, and was developed by the first author in collaboration with Bao Le Hung and Zhixiang Wu. 
It depends on the structure of \(\mathcal O_F/\varpi^2\): when \(F=\mathbb Q_p\), the corresponding filtration is linear; see Remark~\ref{rmk:Qp}. We depict the submodule structure of a representation by the Hasse diagram of its cosocle filtration. The vertices represent composition factors, with successive layers arranged from top to bottom. An edge joining an upper vertex to a lower vertex represents
a non-split extension of the representation at the upper vertex by that at the lower vertex; see \cite{Alp80}.

\begin{proposition}\label{diamond-prop}
For $\sigma$ and $\tau$ as above, the Hasse diagram of $\mathrm{Ind}_I^K V_2(\sigma)^{(\Pi)}$
is
\[
\begin{tikzcd}
   & {\tau^{[s]}}& \\
    {\sigma^{[s]}} & {} & \tau \\
    & {\sigma}
    \arrow[no head, from=1-2, to=2-1]
    \arrow[no head, from=1-2, to=2-3]
    \arrow[no head, from=2-1, to=3-2]
    \arrow[no head, from=2-3, to=3-2]
\end{tikzcd}.
\]
\end{proposition}

\begin{proof}
The constraint $3\leq r\leq p-2$ ensures that all four weights above are
generic and pairwise non-isomorphic. Specializing~\eqref{V2} to
$\rho=\sigma$ and twisting by $\Pi$ gives
\[
0 \longrightarrow \chi^{(\Pi)} \longrightarrow V_2(\sigma)^{(\Pi)}
\longrightarrow \chi^{(\Pi)} \alpha \longrightarrow 0.
\]
Applying $\mathrm{Ind}_I^K$ yields an exact sequence of
$K$-representations
\begin{equation}\label{ses}
0 \longrightarrow \Ind_I^K\chi^{(\Pi)}
\longrightarrow \Ind_I^K V_2(\sigma)^{(\Pi)}
\longrightarrow \Ind_I^K(\chi^{(\Pi)}\alpha)
\longrightarrow 0.
\end{equation}
By the \(f=1\) case of \cite[Theorem~2.4]{BP12}, \(\Ind_I^K\chi^{(\Pi)}\) is a non-split
extension of \(\sigma^{[s]}\) by \(\sigma\), and \(\Ind_I^K(\chi^{(\Pi)}\alpha)\) is a non-split
extension of \(\tau^{[s]}\) by \(\tau\). 

Since \(V_2(\sigma)\) is \(K_1\)-invariant,
\(V_2(\sigma)^{(\Pi)}\) is invariant under \(\Pi^{-1} K_1\Pi\), hence under
\(K_2\subseteq \Pi^{-1} K_1\Pi\). As \(K_2\) is normal in \(K\), this makes
\(\Ind_I^K V_2(\sigma)^{(\Pi)}\) \(K_2\)-invariant. The same is true for
\(\Ind_I^K\chi^{(\Pi)}\) and \(\Ind_I^K(\chi^{(\Pi)}\alpha)\), so \eqref{ses} can be viewed as an exact sequence of
\(K/K_2\)-representations.

\medskip
\noindent\textbf{Step 1.}
\textit{$\mathrm{soc}_{E[K]}(\mathrm{Ind}_I^K V_2(\sigma)^{(\Pi)}) =
\sigma$ and $\mathrm{cosoc}_{E[K]}(\mathrm{Ind}_I^K V_2(\sigma)^{(\Pi)}) =
\tau^{[s]}$.}

\medskip
Indeed, we have
\begin{align*}
\mathrm{soc}_{E[K]}\bigl(\mathrm{Ind}_I^K V_2(\sigma)^{(\Pi)}\bigr)
&\cong \bigoplus_{\rho} \rho \otimes_E \mathrm{Hom}_K\bigl(\rho,\, \mathrm{Ind}_I^K V_2(\sigma)^{(\Pi)}\bigr) \\
&\cong \bigoplus_{\rho} \rho \otimes_E \mathrm{Hom}_I\bigl(\rho|_I,\, V_2(\sigma)^{(\Pi)}\bigr) \\
&\cong \bigoplus_{\rho} \rho \otimes_E \mathrm{Hom}_I\bigl(\rho|_I,\, (V_2(\sigma)^{(\Pi)})^{K_1}\bigr) \\
&\cong \bigoplus_{\rho} \rho\otimes_E \mathrm{Hom}_I\bigl(\rho|_I,\, \chi^{(\Pi)}\bigr) \\
&\cong \bigoplus_{\rho} \rho\otimes_E \mathrm{Hom}_I\bigl(\rho_{I_1},\, \chi^{(\Pi)}\bigr) \\
&\cong \bigoplus_{\rho} \rho\otimes_E \mathrm{Hom}_I\bigl((\rho^{I_1})^{(\Pi)},\, \chi^{(\Pi)}\bigr) \\
&\cong \bigoplus_{\rho} \rho\otimes_E \mathrm{Hom}_I\bigl(\rho^{I_1},\, \chi\bigr) \\
&\cong \sigma,
\end{align*} 
where 
\begin{itemize}
\item the sums range over Serre weights $\rho$; 
\item the first isomorphism comes from evaluation; 
\item the second isomorphism comes from the Frobenius reciprocity; 
\item the third isomorphism comes from the fact that $K_1$ acts trivially on
$\rho$; 
\item the fourth isomorphism comes from $\bigl(V_2(\sigma)^{(\Pi)}\bigr)^{K_1} \cong \chi^{(\Pi)}$; 
\item the fifth isomorphism comes from the fact that $I_1$ acts trivially on $\chi^{(\Pi)}$; and
\item the last three isomorphisms are clear. 
\end{itemize}

The dual argument gives $\mathrm{cosoc}_{E[K]}\bigl(\mathrm{Ind}_I^K
V_2(\sigma)^{(\Pi)}\bigr) = \tau^{[s]}$.

\medskip
\noindent\textbf{Step 2.}
\textit{Writing $W := \mathrm{rad}(\mathrm{Ind}_I^K V_2(\sigma)^{(\Pi)})/\sigma$,
we have $W \cong \sigma^{[s]} \oplus \tau$.}

\medskip
The $K$-representation $W$ fits into the following exact sequence
\[
0 \longrightarrow \sigma^{[s]} \longrightarrow W \longrightarrow \tau \longrightarrow 0,
\]
and hence determines a class in $\mathrm{Ext}^1_{K/K_2}(\tau,
\sigma^{[s]})$. We use the following lemma.

\begin{lemma}\label{ext lem}
Let \(\rho\) and \(\rho'\) be non-isomorphic Serre weights, and
suppose that \[
\Hom_\Gamma\!\left(
\rho,
(\Sym^2E^2\otimes\Det^{-1})\otimes\rho'
\right)=0.
\]
Then the natural inflation map
\(
\Ext^1_\Gamma(\rho,\rho')
\rightarrow
\Ext^1_{K/K_2}(\rho,\rho')
\) is an isomorphism.
\end{lemma}

\begin{proof}
The proof is similar to that of \cite[Corollary~5.2]{BP12}. 
Since
\(K/K_2\cong (K_1/K_2)\rtimes\Gamma\), and \(K_1/K_2\) acts trivially on the
Serre weights \(\rho\) and \(\rho'\), the Lyndon--Hochschild--Serre spectral sequence gives
\[
0\longrightarrow
\Ext^1_\Gamma(\rho,\rho')
\longrightarrow
\Ext^1_{K/K_2}(\rho,\rho')
\longrightarrow
\Hom_\Gamma\bigl(\rho,H^1(K_1/K_2,\rho')\bigr).
\]
Thus, it suffices to prove that the final Hom group vanishes.

Under the identification
\(K_1/K_2=1+\varepsilon M_2(\mathbb F_p)\), the group \(K_1/K_2\) is the
additive group of \(M_2(\mathbb F_p)\), with \(\Gamma\) acting by
conjugation. Hence
\[
H^1(K_1/K_2,\rho')
\cong
\Hom_{\mathbb F_p}(K_1/K_2,\rho')
\cong
M_2(E)^\vee\otimes_E\rho'.
\]
The trace pairing identifies \(M_2(E)^\vee\) with \(M_2(E)\), and
\[
M_2(E)\cong E\oplus(\Sym^2E^2\otimes\Det^{-1})
\]
as a \(\Gamma\)-representation. Therefore, 
\[
\Hom_\Gamma\bigl(\rho,H^1(K_1/K_2,\rho')\bigr)
\cong
\Hom_\Gamma(\rho,\rho')
\oplus
\Hom_\Gamma\!\left(
\rho,
(\Sym^2E^2\otimes\Det^{-1})\otimes\rho'
\right).
\]
The two summands vanish by the hypotheses. Hence, the inflation map
is an isomorphism.
\end{proof}

We apply Lemma~\ref{ext lem} with \(\rho=\tau\) and $\rho'=\sigma^{[s]}$. The weights \(\tau\) and \(\sigma^{[s]}\) are non-isomorphic. It
remains to verify that $\Hom_\Gamma\!\left(
\tau,
(\Sym^2E^2\otimes\det^{-1})\otimes\sigma^{[s]}
\right)=0.
$ The translation principle (or Clebsch--Gordan) gives
\begin{align*}
(\mathrm{Sym}^2 E^2 \otimes \det\nolimits^{-1}) \otimes \sigma^{[s]}
&\cong (\mathrm{Sym}^2 E^2 \otimes \det\nolimits^{-1})
 \otimes (\mathrm{Sym}^{p-1-r} E^2 \otimes \det\nolimits^a) \\
&\cong (\mathrm{Sym}^2 E^2 \otimes \mathrm{Sym}^{p-1-r} E^2)
 \otimes \det\nolimits^{a-1} \\
&\cong
\bigl(\mathrm{Sym}^{p+1-r} E^2 \otimes \det\nolimits^{a-1}\bigr)
\oplus
\bigl(\mathrm{Sym}^{p-1-r} E^2 \otimes \det\nolimits^{a}\bigr) \\
&\qquad{}\oplus
\bigl(\mathrm{Sym}^{p-3-r} E^2 \otimes \Det^{a+1}\bigr) \\
&=
\tau^{[s]} \oplus \sigma^{[s]} \oplus F(p-2+b,a+1),
\end{align*}
where the last summand is omitted when \(p-2+b<a+1\), equivalently when
\(r=p-2\). We already know that
\(\tau\) is not isomorphic to \(\tau^{[s]}\) or \(\sigma^{[s]}\), and it is
not isomorphic to \(F(p-2+b,a+1)\) since the determinant exponents differ by
\((a+1)-(b+1)=r\), which is not divisible by \(p-1\). Hence, the required Hom
group vanishes, and thus Lemma~\ref{ext lem} gives an isomorphism
\[
\Ext^1_\Gamma(\tau,\sigma^{[s]})
\xrightarrow{\sim}
\Ext^1_{K/K_2}(\tau,\sigma^{[s]}).
\]
Consequently, the extension $W$ is the inflation of the
$\Gamma$-representation $W|_\Gamma$, so it suffices to show $W|_\Gamma$
splits.

Let $B \subseteq \Gamma$ be the upper triangular Borel and $U$ its unipotent
radical. Identifying $\Gamma$ with its image in $K/K_2$ under the section (\ref{section}), we have $K/K_2 = (I/K_2)\cdot\Gamma$, and $I/K_2 \cap \Gamma = B$. Using these two equalities, we get that for any $I/K_2$-representation $\eta$,
\begin{equation}\label{ind-res}
\bigl(\Ind_I^K\eta\bigr)\big|_\Gamma
\cong
\bigl(\Ind_{I/K_2}^{K/K_2}\eta\bigr)\big|_\Gamma
\cong
\Ind_B^\Gamma(\eta|_B).
\end{equation}
Restricting~\eqref{ses} to $\Gamma$ and applying
\eqref{ind-res} to each term gives
\begin{equation}\label{ses2}
    0 \to \mathrm{Ind}_B^\Gamma (\chi^{(\Pi)}|_B) \to \mathrm{Ind}_B^\Gamma (V_2(\sigma)^{(\Pi)}|_B) \to \mathrm{Ind}_B^\Gamma(\chi^{(\Pi)}\alpha|_B) \to 0.
\end{equation}
We claim $V_2(\sigma)^{(\Pi)}|_B$ is semisimple. Assuming the claim, we get that (\ref{ses2}) splits and thus $W|_\Gamma$ splits.

It remains to check that $V_2(\sigma)^{(\Pi)}|_B$ is semisimple. Since $V_2(\sigma)$ is
$K_1$-invariant, $V_2(\sigma)^{(\Pi)}$ is invariant under $\Pi^{-1} K_1\Pi$, and a direct computation shows $U \subseteq \Pi^{-1} K_1\Pi/K_2$,
so $U$ acts trivially on $V_2(\sigma)^{(\Pi)}|_B$. Thus, $V_2(\sigma)^{(\Pi)}|_B$ is a $B/U$-representation. Since $|B/U| = (p-1)^2$ is coprime to $p$, it is semisimple.

\medskip
Combining Steps~1 and~2 with $\mathrm{JH}(\mathrm{Ind}_I^K V_2(\sigma)^{(\Pi)}) = \{\sigma,\,
\sigma^{[s]},\, \tau,\, \tau^{[s]}\}$, we obtain the diamond-shaped
Hasse diagram asserted in the proposition.
\end{proof}

\begin{remark}\label{rmk:Qp}
The hypothesis \(F\neq\mathbb Q_p\) in Proposition~\ref{diamond-prop} is
essential. If \(F=\mathbb Q_p\), then the diamond is replaced by the following linear Hasse diagram
\[
\begin{matrix}
\tau^{[s]}\\
\vert\\
\tau\\
\vert\\
\sigma^{[s]}\\
\vert\\
\sigma.
\end{matrix}
\]
Indeed, if the middle extension split, its \(\tau\)-summand would lift to a
non-split extension of \(\tau\) by \(\sigma\). This contradicts
\[
\Ext^1_K(\tau,\sigma)
\cong
\Ext^1_\Gamma(\tau,\sigma)
=
0
\]
for \(F=\mathbb Q_p\); see \cite[Lemma~2.10(ii)]{HuWang} and
\cite[Corollary~5.6(i)]{BP12}.
\end{remark}

\section{Construction of the diagram}

To construct the desired diagram, we will use
Proposition~\ref{diamond-prop} three more times.
We henceforth impose the strengthened hypothesis
\begin{equation}\label{hypothesis}
5\leq r\leq p-2.
\end{equation}
Let \(\kappa\) denote the unique Serre weight satisfying
\(\chi(\kappa)=\chi\alpha^{-2}\). Explicitly,
\(\kappa=F(a-2,b+2)\), and hence
\(\kappa^{[s]}=F(p+1+b,a-2)\). We will apply Proposition~\ref{diamond-prop} with
\(\sigma\) replaced, respectively, by
\[
\tau^{[s]}=F(p+b,a-1),\quad
\tau=F(a-1,b+1),\quad \text{and} \quad
\kappa^{[s]}=F(p+1+b,a-2).
\]
Their coordinate differences are \(p+1-r\), \(r-2\),
and \(p+3-r\), respectively, all of which lie between
\(3\) and \(p-2\) by \eqref{hypothesis}. The proposition therefore applies to all three weights
and shows that the Hasse diagrams for the
\(K\)-representations
\[
\Ind_I^K V_2(\tau^{[s]})^{(\Pi)},\quad
\Ind_I^K V_2(\tau)^{(\Pi)},\quad \text{and}\quad
\Ind_I^K V_2(\kappa^{[s]})^{(\Pi)}
\]
are, respectively,
\[
\begin{tikzcd}
   & \sigma & \\
    {\sigma^{[s]}} & {} & \tau \\
    & {\tau^{[s]}}
    \arrow[no head, from=1-2, to=2-1]
    \arrow[no head, from=1-2, to=2-3]
    \arrow[no head, from=2-1, to=3-2]
    \arrow[no head, from=2-3, to=3-2]
\end{tikzcd},
\quad
\begin{tikzcd}
    & {\kappa^{[s]}} & \\
    {\tau^{[s]}} & {} & \kappa \\
    & \tau
    \arrow[no head, from=1-2, to=2-1]
    \arrow[no head, from=1-2, to=2-3]
    \arrow[no head, from=2-1, to=3-2]
    \arrow[no head, from=2-3, to=3-2]
\end{tikzcd},
\quad \text{and}\quad
\begin{tikzcd}
    & \tau & \\
    {\tau^{[s]}} & {} & \kappa \\
    & {\kappa^{[s]}}
    \arrow[no head, from=1-2, to=2-1]
    \arrow[no head, from=1-2, to=2-3]
    \arrow[no head, from=2-1, to=3-2]
    \arrow[no head, from=2-3, to=3-2]
\end{tikzcd}.
\]

Let $D_0^\#$ be the $K$-representation defined by
\[
D_0^\# :=
\sigma \;\oplus\;
\frac{\mathrm{Ind}_I^K V_2(\sigma)^{(\Pi)}}{\sigma}
\;\oplus\;
\frac{\mathrm{Ind}_I^K V_2(\tau^{[s]})^{(\Pi)}}{\tau^{[s]}}
\;\oplus\;
\frac{\mathrm{Ind}_I^K V_2(\tau)^{(\Pi)}}{\tau}
\;\oplus\;
\frac{\mathrm{Ind}_I^K V_2(\kappa^{[s]})^{(\Pi)}}{\kappa^{[s]}}
\;\oplus\;
\kappa^{[s]}.
\]
Fix \(K\)-isomorphisms
\begin{equation}\label{socle-isomorphisms}
\begin{aligned}
\sigma^{[s]}\oplus\tau
&\xrightarrow{\sim}
\soc_{E[K]}\!\left(\frac{\Ind_I^K V_2(\sigma)^{(\Pi)}}{\sigma}\right),\\
\sigma^{[s]}\oplus\tau
&\xrightarrow{\sim}
\soc_{E[K]}\!\left(\frac{\Ind_I^K V_2(\tau^{[s]})^{(\Pi)}}{\tau^{[s]}}\right),\\
\tau^{[s]}\oplus\kappa
&\xrightarrow{\sim}
\soc_{E[K]}\!\left(\frac{\Ind_I^K V_2(\tau)^{(\Pi)}}{\tau}\right),\,\text{and}\\
\tau^{[s]}\oplus\kappa
&\xrightarrow{\sim}
\soc_{E[K]}\!\left(\frac{\Ind_I^K V_2(\kappa^{[s]})^{(\Pi)}}{\kappa^{[s]}}\right).
\end{aligned}
\end{equation}
Denote the images under the first two isomorphisms in
\eqref{socle-isomorphisms} by
\(\sigma^{[s]}_a\oplus\tau_a\) and
\(\sigma^{[s]}_b\oplus\tau_b\), respectively, and those under
the last two by \(\tau^{[s]}_b\oplus\kappa_b\) and
\(\tau^{[s]}_a\oplus\kappa_a\), respectively. Thus the subscripts \(a\) and \(b\) distinguish the distinct isomorphic submodules of \(D_0^\#\). For \(\rho\in\{\sigma^{[s]},\tau,\tau^{[s]},\kappa\}\)
and \(v\in\rho\), write \(v_a\) and \(v_b\) for its images in \(\rho_a\) and \(\rho_b\), respectively, under the fixed isomorphisms above. It is convenient to understand $D_0^\#$ via the Hasse diagram of each summand
\[
\sigma \;\oplus\;
\begin{tikzcd}
	{\tau^{[s]}} & \\
	{\sigma^{[s]}_a} & {\tau_a}
	\arrow[no head, from=1-1, to=2-1]
	\arrow[no head, from=1-1, to=2-2]
\end{tikzcd}
\;\oplus\;
\begin{tikzcd}
	& \sigma \\
	{\sigma^{[s]}_b} & {\tau^{}_{b}}
	\arrow[no head, from=1-2, to=2-1]
	\arrow[no head, from=1-2, to=2-2]
\end{tikzcd}
\;\oplus\;
\begin{tikzcd}
	{\kappa^{[s]}} & \\
	{\tau^{[s]}_b} & {\kappa_b}
	\arrow[no head, from=1-1, to=2-1]
	\arrow[no head, from=1-1, to=2-2]
\end{tikzcd}
\;\oplus\;
\begin{tikzcd}
	& \tau \\
	{\tau^{[s]}_a} & {\kappa_a}
	\arrow[no head, from=1-2, to=2-1]
	\arrow[no head, from=1-2, to=2-2]
\end{tikzcd}
\;\oplus\;
\kappa^{[s]}.
\]
By abuse of notation, we will sometimes denote a submodule of
\(D_0^\#\) by its corresponding Hasse subdiagram.

We also define, up to isomorphism, an \(I\)-representation
\(D_1\) over \(E\) by
\begin{equation*}
D_1\cong V_2(\sigma)\oplus V_2(\sigma)^{(\Pi)}
\oplus\mathcal C\oplus \mathcal C^{(\Pi)} \oplus V_2(\kappa^{[s]})^{(\Pi)}
\oplus V_2(\kappa^{[s]}),
\end{equation*}
where the \(I\)-representation \(\mathcal C\) is defined as follows.
Choose \(I\)-equivariant injections $j_1:\chi\alpha^{-1}\hookrightarrow V_2(\tau)$ and $j_2:\chi\alpha^{-1}\hookrightarrow V_2(\tau^{[s]})^{(\Pi)}$, and set
\[
\mathcal C:=
\operatorname{coker}\!\left(
\chi\alpha^{-1}
\xrightarrow{(j_1,j_2)}
V_2(\tau)\oplus V_2(\tau^{[s]})^{(\Pi)}
\right).
\]
The Hasse diagram for $D_1$ is given by \[
\begin{tikzcd}[column sep=0.35em, row sep=1.5em]
{\chi \alpha^{-1}} \arrow[d, no head] & \arrow[d, "\oplus" description, draw=none] & {\chi^{(\Pi)}\alpha} \arrow[d, no head] & \arrow[d, "\oplus" description, draw=none] & \chi \arrow[dr, no head] & {\chi\alpha^{-2}} \arrow[d, no head] & \arrow[d, "\oplus" description, draw=none] & {\chi^{(\Pi)}} \arrow[d, no head] & {\chi^{(\Pi)}\alpha^2} \arrow[dl, no head] & \arrow[d, "\oplus" description, draw=none] & {\chi\alpha^{-1}} \arrow[d, no head] & \arrow[d, "\oplus" description, draw=none] & {\chi^{(\Pi)}\alpha} \arrow[d, no head] \\
\chi & {} & {\chi^{(\Pi)}} & {} & {} & {\chi\alpha^{-1}} & {} & {\chi^{(\Pi)}\alpha} & {} & {} & {\chi\alpha^{-2}} & {} & {\chi^{(\Pi)}\alpha^2}
\end{tikzcd}.
\]

\begin{lemma}
$D_0^\#|_I$ contains an $I$-subrepresentation isomorphic to $D_1$.
\end{lemma}

\begin{proof}
By definition, the first and last summands of $D_1$ are isomorphic to $V_2(\sigma)$ and $V_2(\kappa^{[s]})$, which inject into the first and last summands of $D_0^\#|_I$ -- namely $\sigma|_I$ and $\kappa^{[s]}|_I$. For the remaining four summands we use the following claim.

\medskip
\noindent\textbf{Claim.} For $\rho \in \{\sigma, \tau^{[s]}, \tau, \kappa^{[s]}\}$, the composition
\[
\varphi_\rho : V_2(\rho)^{(\Pi)} \hookrightarrow
\mathrm{Ind}_I^K V_2(\rho)^{(\Pi)}
\twoheadrightarrow
\frac{\mathrm{Ind}_I^K V_2(\rho)^{(\Pi)}}{\rho}
\]
of the natural inclusion with the quotient map is injective.

\medskip
\noindent\textit{Proof of Claim.}
In the commutative diagram
\[
\begin{tikzcd}
V_2(\rho)^{(\Pi)} \arrow[hook]{r}
& \mathrm{Ind}_I^K V_2(\rho)^{(\Pi)} \arrow[two heads]{r}
& \dfrac{\mathrm{Ind}_I^K V_2(\rho)^{(\Pi)}}{\rho} \\
\chi(\rho)^{(\Pi)} \arrow[hook]{u} \arrow[hook]{r}
& \mathrm{Ind}_I^K \chi(\rho)^{(\Pi)} \arrow[hook]{u} \arrow[two heads]{r}
& \rho^{[s]} \arrow[hook]{u}
\end{tikzcd}
\]
the vertical maps are inclusions, and the bottom row is the restriction of the top. 
Since $\mathrm{Ind}_I^K \chi(\rho)^{(\Pi)} \twoheadrightarrow \rho^{[s]}$ is nonzero, Frobenius reciprocity makes the composite 
\[\chi(\rho)^{(\Pi)} \hookrightarrow \mathrm{Ind}_I^K\chi(\rho)^{(\Pi)} \twoheadrightarrow \rho^{[s]}\]
nonzero, hence injective. As $\mathrm{soc}_{E[I]}(V_2(\rho)^{(\Pi)}) = \chi(\rho)^{(\Pi)}$, this forces $\varphi_\rho$ to be injective.

\medskip
The claim also shows that each length-three summand of
\(D_0^\#\) has nontrivial \(K_1\)-action, since
\(V_2(\rho)^{(\Pi)}\) does. The \(K_1\)-invariants of each such
summand therefore form a proper \(K\)-submodule containing its
socle. Since this socle is the unique maximal submodule, the
\(K_1\)-invariants equal the socle. As the simple summands of \(D_0^\#\) are \(K_1\)-trivial, we obtain $(D_0^\#)^{K_1}=\operatorname{soc}_{E[K]}D_0^\#$ and $(D_0^\#)^{I_1}=\bigl(\operatorname{soc}_{E[K]}D_0^\#\bigr)^{I_1}$.

Applying the claim with \(\rho=\sigma\) and
\(\rho=\kappa^{[s]}\) gives the desired injections of
\(V_2(\sigma)^{(\Pi)}\) and \(V_2(\kappa^{[s]})^{(\Pi)}\)
into \(\mathrm{Ind}_I^K V_2(\sigma)^{(\Pi)}/\sigma\) and
\(\mathrm{Ind}_I^K V_2(\kappa^{[s]})^{(\Pi)}/\kappa^{[s]}\),
respectively.

For the third summand of \(D_1\), taking
\(\rho=\tau^{[s]}\) in the claim gives an injection of
\(V_2(\tau^{[s]})^{(\Pi)}\) into
\(\mathrm{Ind}_I^K V_2(\tau^{[s]})^{(\Pi)}/\tau^{[s]}\).
 The same target also contains a copy of $\chi\alpha^{-1}$ via the inclusion of $V_2(\tau)$ into $\tau$, which injects into $\mathrm{Ind}_I^K V_2(\tau^{[s]})^{(\Pi)}/\tau^{[s]}$. Since $\chi\alpha^{-1}$ has multiplicity one in the $E[I]$-socle of this target, these two copies of $\chi\alpha^{-1}$ must coincide, and we get an injection of $\mathcal C$ into $\mathrm{Ind}_I^K V_2(\tau^{[s]})^{(\Pi)}/\tau^{[s]}$. An identical argument, with the roles of $\tau$ and $\tau^{[s]}$ exchanged and $\chi\alpha^{-1}$ replaced by $\chi^{(\Pi)}\alpha$, gives an injection of $\mathcal C^{(\Pi)}$ into $\mathrm{Ind}_I^K V_2(\tau)^{(\Pi)}/\tau$. Summing all six injections yields the desired $I$-linear embedding $D_1 \hookrightarrow D_0^\#|_I$.
\end{proof}
We henceforth identify $D_1$ with its image in $D_0^\#|_I$ under this embedding.

For each \(i\in\mathbb{Z}\), let \(D_0^\#(i)\) be a copy of
\(D_0^\#\), and let $D_1(i)\subseteq D_0^\#(i)|_I
$ be the corresponding copy of \(D_1\subseteq D_0^\#|_I\). Set
\[
D_0^\#(\infty):=\bigoplus_{i\in\mathbb{Z}}D_0^\#(i)
\quad \text{and} \quad
D_1(\infty):=\bigoplus_{i\in\mathbb{Z}}D_1(i).
\]

Now, let \(\widetilde E\) be a division algebra over \(E\) with center
\(\widetilde C\). For any \(E\)-representation \(V\), write
\(\widetilde V:=V\otimes_E\widetilde E\). We also use this convention for the indexed copies and their direct sums. These representations
carry natural commuting left and right \(\widetilde E\)-actions; we
distinguish them by writing scalars on the left and on the right,
respectively.

By abuse of notation, we identify \(D_0^\#\), \(D_1\), \(\sigma\),
etc.\ with their images under the natural injections $D_0^\#\hookrightarrow\widetilde D_0^\#$, $D_1\hookrightarrow\widetilde D_1$, $\sigma\hookrightarrow\widetilde\sigma$, and so on. We make the same identifications for the indexed copies and their direct sums.

Fix a \(\widetilde C\)-valued character \(\xi\) of
\(
\begin{psmallmatrix}
\varpi&0\\
0&\varpi
\end{psmallmatrix}^{\mathbb{Z}}\cong \mathbb{Z}.
\)
We extend the induced \(K\)-actions on
\(\widetilde D_0^\#\) and \(\widetilde D_0^\#(\infty)\), and the
induced \(I\)-actions on \(\widetilde D_1\) and
\(\widetilde D_1(\infty)\), to actions of
\(
KZ\cong
K\times
\begin{psmallmatrix}
\varpi&0\\
0&\varpi
\end{psmallmatrix}^{\mathbb{Z}}
\) and
\(
IZ\cong
I\times
\begin{psmallmatrix}
\varpi&0\\
0&\varpi
\end{psmallmatrix}^{\mathbb{Z}},
\) respectively, by letting
\(
\begin{psmallmatrix}
\varpi&0\\
0&\varpi
\end{psmallmatrix}^{\mathbb{Z}}
\) act through \(\xi\). For brevity, we write
\(
\xi(\varpi):=
\xi\left(
\begin{psmallmatrix}
\varpi&0\\
0&\varpi
\end{psmallmatrix}
\right).
\)

For each \(i\in\mathbb{Z}\), the fixed identification of
\(D_0^\#(i)\) with \(D_0^\#\) induces an isomorphism
\(
\widetilde D_0^\#
\xrightarrow{\sim}
\widetilde D_0^\#(i)
\) of \(KZ\)-representations commuting with the left and right
\(\widetilde E\)-actions and carrying \(\widetilde D_1\) onto
\(\widetilde D_1(i)\). Let
\[
\iota_i:
\widetilde D_0^\#
\xrightarrow{\sim}
\widetilde D_0^\#(i)
\hookrightarrow
\widetilde D_0^\#(\infty)
\]
denote the resulting inclusion of the \(i\)-th copy, and write
\(v_i:=\iota_i(v)\) for \(v\in\widetilde D_0^\#\). Each \(\iota_i\)
commutes with the left and right \(\widetilde E\)-actions. We use the
same \(a\)- and \(b\)-subscripts after extension of scalars to
\(\widetilde E\). Thus, for $\rho\in\{\sigma^{[s]},\tau,\tau^{[s]},\kappa\},
$ \(\widetilde{\rho}_a\) and \(\widetilde{\rho}_b\) denote the corresponding images of \(\widetilde{\rho}\) in \(\widetilde{D}_0^\#\). For \(i,j\in\mathbb Z\) and $\rho\in\{\sigma^{[s]},\tau^{[s]},\kappa\}$, define \(\Delta_{i,j}(\widetilde{\rho})\) to be the image of the
following composition
\[
\widetilde{\rho}
\xrightarrow{\ v\mapsto (v_a,v_b)\ }
\widetilde{\rho}_a\oplus\widetilde{\rho}_b
\xrightarrow{\ \iota_i\oplus\iota_j\ }
\widetilde{D}_0^\#(i)\oplus\widetilde{D}_0^\#(j)
\xrightarrow{\ (x,y)\mapsto x-y\ }
\M .
\]
Equivalently,
\begin{equation}\label{Del}
\Delta_{i,j}(\widetilde{\rho})
=
\{\, (v_a)_i-(v_b)_j \mid v\in\widetilde{\rho}\,\}
\subseteq \M .
\end{equation}

\noindent
Similarly, for \(\lambda_i\in\widetilde{E}^{\times}\), define
\(\Delta_{i,i+1}^{\lambda_i}(\widetilde{\tau})\) to be the image of the following composition
\[
\widetilde{\tau}
\xrightarrow{\ v\mapsto (v_a, v_b\lambda_i)\ }
\widetilde{\tau}_a\oplus\widetilde{\tau}_b
\xrightarrow{\ \iota_i\oplus\iota_{i+1}\ }
\widetilde{D}_0^\#(i)\oplus\widetilde{D}_0^\#(i+1)
\xrightarrow{\ (x,y)\mapsto x-y\ }
\M .
\]
Equivalently,
\begin{equation}\label{DelLam}
\Delta_{i,i+1}^{\lambda_i}(\widetilde{\tau})
=
\{\, (v_a)_i-(v_b)_{i+1}\lambda_i\mid v\in\widetilde{\tau}\,\}
\subseteq \M .
\end{equation}

Fix $\lambda:=(\lambda_i)_{i\in\mathbb Z}
\in\prod_{i\in\mathbb Z}\widetilde{E}^{\times}
$ and define \(R\) to be the \(KZ\)-subrepresentation of \(\M\) given by
\[
R:=
\bigoplus_{i\in\mathbb Z}
\Bigl(
\Delta_{i,i}(\widetilde{\sigma}^{[s]})
\oplus
\Delta_{i,i+1}^{\lambda_i}(\widetilde{\tau})
\oplus
\Delta_{i,i-1}(\widetilde{\tau}^{[s]})
\oplus
\Delta_{i,i}(\widetilde{\kappa})
\Bigr).
\]
Finally, set $\widetilde{D}_0(\infty):=\M/R.
$

\begin{lemma}\label{D1R0}
The intersection of $\widetilde{D}_1(\infty)$ and $R$ inside $\widetilde{D}_0^\#(\infty)$ is zero. 
In particular, the natural map $\phi :\widetilde{D}_1(\infty) \rightarrow \widetilde{D}_0(\infty)
$ is injective.
\end{lemma}

\begin{proof}
To show that $\widetilde{D}_1(\infty) \cap R = 0$ inside
$\widetilde{D}_0^\#(\infty)$, it suffices to show that $\bigl(\widetilde{D}_1(\infty)
\cap R\bigr)^{I_1} = 0$. Indeed, every nonzero smooth representation of the
pro-\(p\) group \(I_1\) over a field of characteristic
\(p\) has nonzero \(I_1\)-invariants.

Since $\bigl(\widetilde{D}_1(\infty) \cap R\bigr)^{I_1} =
\widetilde{D}_1(\infty)^{I_1} \cap R^{I_1}$, with
\[
\widetilde{D}_1(\infty)^{I_1} = \bigoplus_{i \in \mathbb{Z}}\iota_i
\Bigl( \widetilde{\sigma}^{I_1} \oplus (\widetilde{\sigma}^{[s]}_a)^{I_1}
\oplus (\widetilde{\tau}_b)^{I_1} \oplus (\widetilde{\tau}^{[s]}_b)^{I_1}
\oplus (\widetilde{\kappa}_a)^{I_1} \oplus (\widetilde{\kappa}^{[s]})^{I_1} \Bigr),
\]
it suffices to show that the composition
\begin{equation*}
    R^{I_1}\hookrightarrow \widetilde{D}_0^\#(\infty)^{I_1}
\twoheadrightarrow \widetilde{D}_0^\#(\infty)^{I_1}\big/\widetilde{D}_1(\infty)^{I_1} 
\cong \bigoplus_{i \in \mathbb{Z}} \iota_i\Bigl( (\widetilde{\sigma}^{[s]}_b)^{I_1}
\oplus (\widetilde{\tau}_a)^{I_1} \oplus (\widetilde{\tau}^{[s]}_a)^{I_1}
\oplus (\widetilde{\kappa}_b)^{I_1} \Bigr)
\end{equation*}
is injective, where the first map is the inclusion, the second is the
quotient map, and the isomorphism is induced from the projection map. 

Distinct summands of \(R^{I_1}\) map into distinct
summands on the right, so it suffices to check that
the restriction to each summand is injective. A nonzero element of
$\Delta_{i,i}\bigl((\widetilde{\sigma}^{[s]})^{I_1}\bigr)$ maps to a
nonzero element of the $i$-th copy of $(\widetilde{\sigma}^{[s]}_b)^{I_1}$,
and likewise a nonzero element of $\Delta_{i,i}\bigl((\widetilde{\kappa})^{I_1}\bigr)$
maps to a nonzero element of the $i$-th copy of $(\widetilde{\kappa}_b)^{I_1}$.
A nonzero element of $\Delta_{i,i+1}^{\lambda_i}\bigl((\widetilde{\tau})^{I_1}\bigr)$
maps to a nonzero element of the $i$-th copy of $(\widetilde{\tau}_a)^{I_1}$,
and likewise a nonzero element of
$\Delta_{i,i-1}\bigl((\widetilde{\tau}^{[s]})^{I_1}\bigr)$ maps to a nonzero
element of the $i$-th copy of $(\widetilde{\tau}^{[s]}_a)^{I_1}$. This
proves injectivity, and hence the lemma.
\end{proof}

We now extend the \(IZ\)-action on \(\widetilde D_1(\infty)\) to an
\(N_G(I)\)-action. Let \(A\) be
the direct sum of the first, third, and fifth summands of \(D_1\), and let
\(B\) be the direct sum of the second, fourth, and sixth summands. Thus,
\(D_1=A\oplus B\), with \(A\cong B^{(\Pi)}\) as \(E[I]\)-representations. Fix an \(E[I]\)-isomorphism
\(\psi:A\xrightarrow{\sim}B^{(\Pi)}\) given by the direct sum of isomorphisms between the three corresponding summands. Extending scalars gives an
\(E[I]\)-isomorphism
\begin{equation}\label{ABpsi}
\widetilde\psi:=\psi\otimes_E \mathrm{id}_{\widetilde E}:
\widetilde A\xrightarrow{\sim}\widetilde B^{(\Pi)}.
\end{equation}
Since \(\widetilde\psi\) acts as the identity on the second tensor
factor, the map \(\widetilde\psi\) commutes
with the left and right \(\widetilde E\)-actions. There is then a
unique extension of the given \(IZ\)-action on \(\widetilde D_1\) to an
\(N_G(I)\)-action such that the map $v \mapsto (\Pi v)^{(\Pi)}$ for $v \in \widetilde{A}$ is $\widetilde{\psi}$. Moreover, this \(N_G(I)\)-action commutes with
the left and right \(\widetilde E\)-actions.

We define the \(N_G(I)\)-action on \(\widetilde D_1(\infty)\) componentwise
from the action on \(\widetilde D_1\). Explicitly, for
\(v\in\widetilde D_1\) and \(i\in\mathbb Z\), set
\begin{equation}\label{Pi action}
\Pi v_i := (\Pi v)_i .
 \end{equation}
This uniquely determines a smooth \(N_G(I)\)-action on
\(\widetilde D_1(\infty)\) extending the given \(IZ\)-action. Together with the map \(\phi\) from Lemma~\ref{D1R0},
this gives a diagram
\((\widetilde D_0(\infty),\widetilde D_1(\infty),\phi)\)
in the sense of \cite[\S9]{BP12}. A significant point of departure from
much of the previous literature is that this diagram is neither
a basic \(0\)-diagram nor \(K_1\)-invariant. The following proposition, the
analog of \cite[Proposition~2.1]{Le24}, is the key irreducibility statement
for our diagram.

In what follows, all products of \(\lambda_i\) are ordered from
left to right according to increasing indices. Thus, for
\(i\in\mathbb Z\) and \(m\geq 1\), we set \(
\prod_{j=0}^{m-1}\lambda_{i+j}
:=
\lambda_i\lambda_{i+1}\cdots\lambda_{i+m-1}.
\)

\begin{proposition}\label{irred prop}
Let $W^\# \subseteq \widetilde{D}_0^\#(\infty)$ be an $E[KZ]$-subrepresentation such that $R \subsetneq W^\#$ and $W^\# \cap \widetilde{D}_1(\infty)\subseteq \widetilde{D}_0^\#(\infty)$ is stable under the action of $\Pi$. If $\xi(\varpi)^4 \lambda_0 \in E^{\times}$, $\lambda_i \neq \lambda_0$ for all $i \neq 0$, and the $E$-span of
\[
\left\{
\xi(\varpi)^{4m}
\prod_{j=0}^{m-1} \lambda_{j}
\;\middle|\;
m \ge 1
\right\}
\]
is $\widetilde{E}$, then $W^\# = \widetilde{D}_0^\#(\infty)$.
\end{proposition}
\begin{proof}
    We begin with a series of lemmata.

\begin{lemma}\label{propagation lemma}
Let \(W^\#\) be as in Proposition~\ref{irred prop}. If \((c_i)_{i\in\mathbb Z}\in \bigoplus_{i\in\mathbb Z}\widetilde E\) satisfies \(\left(\sum_i\iota_i c_i\right)(\tau_b)\subseteq W^\#\), then, for every \(m\geq1\),
\[
\left(
\xi(\varpi)^{2m}
\sum_i
\iota_{i+m}
\left(
c_i\prod_{j=0}^{m-1}\lambda_{i+j}
\right)
\right)(\tau_b)
\subseteq W^\#
\]
and
\[
\left(
\xi(\varpi)^{2m}
\sum_i\iota_{i-m}c_i
\right)(\tau_b)
\subseteq W^\#.
\]
\end{lemma}

\begin{proof}
The proof proceeds by iterating the following four operations.

\medskip

\noindent
(1) \emph{Application of $\Pi$.} 
Applying the action of $\Pi$ allows us to pass between different summands in $\widetilde{D}_1(\infty)$ and introduces the scalar $\xi(\varpi)$ appearing in the statement. Note that the action of $\Pi$ on each summand coming from $B$ (the direct sum of the second,
fourth, and sixth summands of \(D_1\)) is determined by the fact that $\Pi^2 = \begin{psmallmatrix}
\varpi&0\\
0&\varpi
\end{psmallmatrix}$, and $\begin{psmallmatrix}
\varpi&0\\
0&\varpi
\end{psmallmatrix}$ acts by $\xi(\varpi)$.

\medskip

\noindent
(2) \emph{Use of the gluing relations.}
Since \(R\subseteq W^\#\), the relation $(v_a)_i-(v_b)_{i+1}\lambda_i
\in
R
$ allows us to pass from an \(a\)-copy of \(\widetilde\tau\) at index
\(i\) to a \(b\)-copy at index \(i+1\), with right multiplication
by \(\lambda_i\). More precisely, for \(c\in\widetilde E\), if \(
\iota_i(c\tau_a)\subseteq W^\#\), then \(\iota_{i+1}(c\lambda_i\tau_b)\subseteq W^\#.
\) We use the analogous unscaled gluing relations for
\(\widetilde\sigma^{[s]}\), \(\widetilde\tau^{[s]}\), and
\(\widetilde\kappa\).

\medskip

\noindent
(3) \emph{Generation of $K$-subrepresentations.}
We will use that the $K$-subrepresentations of
\(\widetilde D_0^\#\) generated by
\(V_2(\sigma)^{(\Pi)}\),
\(V_2(\tau^{[s]})^{(\Pi)}\),
\(V_2(\tau)^{(\Pi)}\), and
\(V_2(\kappa^{[s]})^{(\Pi)}\)
have the following Hasse diagrams, respectively:
\[
\begin{tikzcd}
	{\tau^{[s]}} & \\
	{\sigma^{[s]}_a} & {\tau_a}
	\arrow[no head, from=1-1, to=2-1]
	\arrow[no head, from=1-1, to=2-2]
\end{tikzcd}
\quad
\begin{tikzcd}
	& \sigma \\
	{\sigma^{[s]}_b} & {\tau_b}
	\arrow[no head, from=1-2, to=2-1]
	\arrow[no head, from=1-2, to=2-2]
\end{tikzcd}
\quad
\begin{tikzcd}
	{\kappa^{[s]}} & \\
	{\tau^{[s]}_b} & {\kappa_b}
	\arrow[no head, from=1-1, to=2-1]
	\arrow[no head, from=1-1, to=2-2]
\end{tikzcd}
\quad
\begin{tikzcd}
	& \tau \\
	{\tau^{[s]}_a} & {\kappa_a}
	\arrow[no head, from=1-2, to=2-1]
	\arrow[no head, from=1-2, to=2-2]
\end{tikzcd}.
\]
We also use the fact that for any Serre weight $\rho$, the $K$-subrepresentation generated by the $I_1$-invariant line $\rho^{I_1}$ is $\rho$.

\medskip

\noindent
(4) \emph{Closure under subrepresentations.} 
If \(M\subseteq W^\#\) is an \(E[K]\)-submodule, then every \(E[I]\)-submodule of \(M\) is contained in \(W^\#\).

\medskip
We begin with the assumption
\[
\Bigl( \sum_i \iota_i c_i \Bigr)(\tau_b) \subseteq W^\#.
\]
Applying (4) to $V_2(\tau_b)$, followed by (1) and (3), we obtain
\[
\Bigl( \sum_i \iota_i c_i \Bigr)
\left(
\begin{tikzcd}
\kappa^{[s]} & \\
\tau_b^{[s]} & \kappa_b
\arrow[no head, from=1-1, to=2-1]
\arrow[no head, from=1-1, to=2-2]
\end{tikzcd}
\right)
\subseteq W^\#.
\]
Applying (4) to $V_2(\tau_b^{[s]})$, followed by (1) and (3), we get
\[
\xi(\varpi) \Bigl( \sum_i \iota_i c_i \Bigr)
\left(
\begin{tikzcd}
& \sigma \\
\sigma_b^{[s]} & \tau_b
\arrow[no head, from=1-2, to=2-1]
\arrow[no head, from=1-2, to=2-2]
\end{tikzcd}
\right)
\subseteq W^\#.
\]
Applying (4) to $\sigma_b^{[s]}$ and then (2), we obtain
\[
\xi(\varpi) \Bigl( \sum_i \iota_i c_i \Bigr)(\sigma_a^{[s]}) \subseteq W^\#.
\]
Applying (4) to $(\sigma_a^{[s]})^{I_1}$, followed by (1) and (3), we get
\[
\xi(\varpi)^2 \Bigl( \sum_i \iota_i c_i \Bigr)(\sigma) \subseteq W^\#.
\]
Applying (4) to $V_2(\sigma)$, followed by (1) and (3), we have
\[
\xi(\varpi)^2 \Bigl( \sum_i \iota_i c_i \Bigr)
\left(
\begin{tikzcd}
\tau^{[s]} & \\
\sigma_a^{[s]} & \tau_a
\arrow[no head, from=1-1, to=2-1]
\arrow[no head, from=1-1, to=2-2]
\end{tikzcd}
\right)
\subseteq W^\#.
\]
Finally, applying (4) to $\tau_a$ followed by (2) gives
\[
\xi(\varpi)^2 \Bigl( \sum_i \iota_{i+1} c_i \lambda_i \Bigr)(\tau_b) \subseteq W^\#.
\]
Repeating this argument inductively, we obtain for all $m \ge 1$,
\[
\Bigl( \xi(\varpi)^{2m} \sum_i \iota_{i+m} \bigl(c_i \prod_{j=0}^{m-1} \lambda_{i+j} \bigr) \Bigr)(\tau_b) \subseteq W^\#.
\]

To prove the second containment, again start by applying (4) to $V_2(\tau_b)$ followed by (1) and (3) to get
\[
\Bigl( \sum_i \iota_i c_i \Bigr)
\left(
\begin{tikzcd}
\kappa^{[s]} & \\
\tau_b^{[s]} & \kappa_b
\arrow[no head, from=1-1, to=2-1]
\arrow[no head, from=1-1, to=2-2]
\end{tikzcd}
\right)
\subseteq W^\#.
\]
Applying (4) to $\kappa_b$ followed by (2) gives
\[
\Bigl( \sum_i \iota_i c_i \Bigr)(\kappa_a) \subseteq W^\#.
\]
Applying (4) to $\kappa_a^{I_1}$, followed by (1) and (3), we get
\[
\Bigl( \sum_i \iota_i c_i \Bigr)(\kappa^{[s]}) \subseteq W^\#.
\]
Applying (4) to $V_2(\kappa^{[s]})$, followed by (1) and (3), we obtain
\[
\xi(\varpi) \Bigl( \sum_i \iota_i c_i \Bigr)
\left(
\begin{tikzcd}
& \tau \\
\tau_a^{[s]} & \kappa_a
\arrow[no head, from=1-2, to=2-1]
\arrow[no head, from=1-2, to=2-2]
\end{tikzcd}
\right)
\subseteq W^\#.
\]
Applying (4) to $\tau_a^{[s]}$ followed by (2) gives
\[
\xi(\varpi) \Bigl( \sum_i \iota_{i-1} c_i \Bigr)(\tau_b^{[s]}) \subseteq W^\#.
\]
Applying (4) to $V_2(\tau_b^{[s]})$, followed by (1) and (3), we get
\[
\xi(\varpi)^2 \Bigl( \sum_i \iota_{i-1} c_i \Bigr)
\left(
\begin{tikzcd}
& \sigma \\
\sigma_b^{[s]} & \tau_b
\arrow[no head, from=1-2, to=2-1]
\arrow[no head, from=1-2, to=2-2]
\end{tikzcd}
\right)
\subseteq W^\#.
\]
Finally, applying (4) to $\tau_b$ gives
\[
\xi(\varpi)^2 \Bigl( \sum_i \iota_{i-1} c_i \Bigr)(\tau_b) \subseteq W^\#.
\]
Repeating this argument and proceeding by induction on $m$, we obtain, for all $m \ge 1$,
\[
\Bigl( \xi(\varpi)^{2m} \sum_i \iota_{i-m} c_i \Bigr)(\tau_b) \subseteq W^\#.
\]
\end{proof}

\begin{lemma}\label{socle surjection}
The quotient map
\(
q:\widetilde D_0^\#(\infty)\rightarrow
\widetilde D_0(\infty)
\)
induces a surjection $$\operatorname{soc}_{E[K]}\widetilde D_0^\#(\infty)
\rightarrow
\operatorname{soc}_{E[K]}\widetilde D_0(\infty).
$$ \end{lemma}

\begin{proof}

Let
\(
\rho\subseteq\operatorname{soc}_{E[K]}\widetilde D_0(\infty)
\)
be a Serre weight, and let $f:\rho\hookrightarrow\widetilde D_0(\infty)
$ denote the inclusion. Let
\(
p:P(\rho)\twoheadrightarrow\rho
\)
be a projective cover of $\rho$ in the category of
\(E[K/K_2]\)-representations. By projectivity, the composite
\(f\circ p\) lifts through \(q\) to an \(E[K]\)-map $F:P(\rho)\rightarrow\widetilde D_0^\#(\infty).
$ Since \(q\circ F=f\circ p\),
\begin{equation}\label{F(rad)in R}
F(\operatorname{rad}P(\rho))\subseteq R.
\end{equation}

For each \(i\in\mathbb Z\), let
\(
\operatorname{pr}_i:
\widetilde D_0^\#(\infty)\rightarrow
\widetilde D_0^\#(i)
\)
denote the natural projection. Suppose that
\(\operatorname{pr}_i\circ F\neq0\).
The length-three summand \(Q\) of \(D_0^\#\) with cosocle
\(\rho\), if it exists, is unique. Since \(\rho\) occurs
in \(Q\) with multiplicity one, projectivity gives
\(\dim_E\operatorname{Hom}_{E[K]}(P(\rho),Q)=1\).
As \(P(\rho)\) is finite dimensional, every nonzero map
\(P(\rho)\to\widetilde Q\) therefore has the kernel of a
fixed surjection \(P(\rho)\twoheadrightarrow Q\), which
is contained in \(\ker p\). All other components factor
through \(p\), since their targets modulo their socles
have no composition factor isomorphic to \(\rho\).
Consequently, \(\operatorname{im}(\operatorname{pr}_i\circ F)\)
is isomorphic either to \(\rho\) or to \(Q\).
We show that the latter possibility cannot occur.

Suppose first that \(\rho=\tau^{[s]}\) and that, for some \(i\in\mathbb Z\),
\(\operatorname{im}(\operatorname{pr}_i\circ F)\) is non-simple.
Then the Hasse diagram for $\operatorname{im}(\operatorname{pr}_i\circ F)$ is $\begin{tikzcd}[row sep=0.7em,column sep=1.2em]\tau^{[s]} & \\
\sigma^{[s]} & \tau
\arrow[no head, from=1-1, to=2-1]
\arrow[no head, from=1-1, to=2-2]
\end{tikzcd}$, and consequently $(\operatorname{pr}_i\circ F)
\bigl(\operatorname{rad}P(\tau^{[s]})\bigr)
\cong
\sigma^{[s]}\oplus\tau$. Consider the projection of \(F(P(\tau^{[s]}))\) onto $\iota_i\left(\begin{tikzcd}[row sep=0.7em,column sep=1.2em]
& \sigma \\
\sigma^{[s]}_b & \tau_b
\arrow[no head, from=1-2, to=2-1]
\arrow[no head, from=1-2, to=2-2]
\end{tikzcd}\right)\subseteq \widetilde D_0^\#(i)$. This projection is zero, since
\(\tau^{[s]}\) is not a Jordan--Hölder constituent of this summand. Hence $(\operatorname{pr}_i\circ F)(\operatorname{rad}P(\tau^{[s]}))\cap \iota_i(\widetilde\sigma_a^{[s]}) \neq 0$ inside $\widetilde D_0^\#(i)$. This contradicts \eqref{F(rad)in R}, as the
\(\sigma^{[s]}\)-part of \(R\) is
\(
\bigoplus_{j\in\mathbb Z}
\Delta_{j,j}(\widetilde\sigma^{[s]}).
\)

The case \(\rho=\sigma\) is the same with the \(a\)- and \(b\)-copies interchanged. The cases \(\rho=\kappa^{[s]}\) and \(\rho=\tau\) are identical, using \(\Delta_{i,i}(\widetilde\kappa)\) in place of \(\Delta_{i,i}(\widetilde\sigma^{[s]})\). Finally, \(\sigma^{[s]}\) and \(\kappa\) occur only in the socle of \(D_0^\#\), so no non-simple image is possible for these two weights. Thus, \(\operatorname{pr}_i\circ F\) annihilates \(\operatorname{rad}P(\rho)\) for every $i\in\mathbb Z$. Hence \( F(\operatorname{rad}P(\rho))=0, \) and \(F\) factors through \(p\), giving a nonzero \(E[K]\)-map \( f^\#:\rho\rightarrow \operatorname{soc}_{E[K]}\widetilde D_0^\#(\infty) \) such that \(q\circ f^\#=f\). This proves the surjectivity on socles.
\end{proof}

\begin{lemma}\label{nonzero hom}
Let \(W^\#\) be as in Proposition~\ref{irred prop}. Suppose that $\xi(\varpi)^4\lambda_0 \in E^\times$ and $\lambda_i \neq \lambda_0$ for all $i \neq 0$. Then there exists \(d_0\in\widetilde E^\times\) such that $\iota_0d_0(\tau_b)\subseteq W^\#$.
\end{lemma}

\begin{proof}
Let $q:\widetilde D_0^\#(\infty)\rightarrow \widetilde D_0(\infty)$ be the quotient map. Since \(R\subsetneq W^\#\), $q(W^\#)\neq 0$. Thus, there exists a nonzero \(E[K]\)-map \(f:\rho\hookrightarrow q(W^\#)\) for some Serre weight
\(\rho\). By Lemma~\ref{socle surjection}, \(f\) lifts to a nonzero \(E[K]\)-map $f^\#:\rho\rightarrow
\operatorname{soc}_{E[K]}\widetilde D_0^\#(\infty)$ such that \(q\circ f^\#=f\). Since \(f(\rho)\subseteq q(W^\#)\) and \(\ker q=R\subseteq W^\#\),
we have \(f^\#(\rho)\subseteq W^\#\). As $\rho$ is isomorphic to a simple $E[K]$-subrepresentation of $D_0^\#(\infty)$, we have the following possible cases:

\medskip
\noindent
\emph{Case 1: \(\rho\cong\sigma\) or \(\kappa^{[s]}\).} Then
\[
f^\#(\rho)=
\left(\sum_i\iota_i c_i\right)(\rho)\subseteq W^\#
\]
for some nonzero family
\((c_i)\in\bigoplus_{i\in\mathbb Z}\widetilde E\).

\medskip
\noindent
\emph{Case 2: \(\rho\cong\sigma^{[s]},\tau,\tau^{[s]}\),
or \(\kappa\).}
Using the fixed isomorphisms in \eqref{socle-isomorphisms}, we can write
\[
f^\#(v)=
\sum_i\iota_i(a_i v_a)+
\sum_i\iota_i(b_i v_b)\in W^\#
\qquad\textup{for all } v\in\rho,
\]
where
\((a_i),(b_i)\in\bigoplus_{i\in\mathbb Z}\widetilde E\).
Set
\[
c_i:=
\begin{cases}
b_i+a_i,
  &\text{if }\rho\cong\sigma^{[s]}\text{ or }\kappa,\\
b_i+a_{i-1}\lambda_{i-1},
  &\text{if }\rho\cong\tau,\\
b_i+a_{i+1},
  &\text{if }\rho\cong\tau^{[s]}.
\end{cases}
\]
Then $f^\#(v)-\sum_i\iota_i(c_i v_b)\in R$ for all $v\in\rho$. Since \(R\subseteq W^\#\), it follows that
\[
\left(\sum_i\iota_i c_i\right)(\rho_b)\subseteq W^\#.
\]
Moreover, \((c_i)\neq0\), since otherwise
\(f^\#(\rho)\subseteq R\), contradicting
\(q\circ f^\#=f\neq0\).

\medskip
Starting from the containment obtained in either case,
the argument in the proof of
Lemma~\ref{propagation lemma} yields a family
\((d_i)\in\bigoplus_{i\in\mathbb Z}\widetilde E\)
such that
\begin{equation}\label{di}
\left(\sum_i\iota_i d_i\right)(\tau_b)\subseteq W^\#.
\end{equation}
At each step in the proof of Lemma~\ref{propagation lemma}, the coefficient family is either unchanged or modified by translating the indices or multiplying the coefficients by invertible scalars. Since \((c_i)\neq0\) in both cases above, the resulting family \((d_i)\) is also nonzero.

Suppose that $\#(d_i):=\#\{i\in\mathbb Z\mid d_i\neq0\}$ is minimal among all nonzero families \((d_i)\) satisfying \eqref{di}. By Lemma~\ref{propagation lemma}, we may assume without loss of generality that \(d_0\neq0\).
Applying the first containment of
Lemma~\ref{propagation lemma} gives
\[
\Bigl(\xi(\varpi)^2
\sum_i\iota_{i+1}d_i\lambda_i\Bigr)(\tau_b)
\subseteq W^\#,
\]
and applying the second containment of
Lemma~\ref{propagation lemma} to shift back by one
index gives
\[
\Bigl(\xi(\varpi)^4
\sum_i\iota_i d_i\lambda_i\Bigr)(\tau_b)
\subseteq W^\#.
\]
Using \(\xi(\varpi)^4\lambda_0\in E^\times\)
(and hence \(\lambda_0\in\widetilde C^\times\)), we obtain
\[
\left(
\xi(\varpi)^4
\sum_i\iota_i d_i(\lambda_i-\lambda_0)
\right)(\tau_b)\subseteq W^\#.
\]
Since \(\lambda_i-\lambda_0\neq0\) for all \(i\neq0\), $\#\bigl(d_i(\lambda_i-\lambda_0)\bigr)=\#(d_i)-1$. Minimality implies \(d_i=0\) for all \(i\neq0\), and hence $\iota_0d_0(\tau_b)\subseteq W^\#$ for some $d_0\in\widetilde E^\times$.
\end{proof}

\begin{lemma}\label{i_0 tilde}
Let \(W^\#\) and $\lambda=(\lambda_i)_{i\in\mathbb Z}
\in\prod_{i\in\mathbb Z}\widetilde{E}^{\times}$ be as in Proposition~\ref{irred prop}. Then \(\iota_0(\widetilde\tau_b)\subseteq W^\#\).
\end{lemma}

\begin{proof}
By Lemma~\ref{nonzero hom}, there exists \(d_0\in\widetilde E^\times\)
such that \(\iota_0d_0(\tau_b)\subseteq W^\#\). By
Lemma~\ref{propagation lemma}, we also have
\[
\iota_0\!\left(\xi(\varpi)^{4m}d_0\prod_{j=0}^{m-1}\lambda_j\right)
(\tau_b)\subseteq W^\#
\]
for all \(m\ge1\). Since \(\xi(\varpi)\in\widetilde C^\times\), we have
\(\xi(\varpi)^{4m}d_0\prod_{j=0}^{m-1}\lambda_j
=d_0 \left( \xi(\varpi)^{4m}\prod_{j=0}^{m-1}\lambda_j\right)\) for all \(m\ge1\). By the choice of \(\lambda\), the \(E\)-span of $\left\{
\xi(\varpi)^{4m}\prod_{j=0}^{m-1}\lambda_j
\;\middle|\;
m\ge1
\right\}$ is \(\widetilde E\). Since left multiplication by \(d_0\in\widetilde E^\times\) is an \(E\)-linear bijection on \(\widetilde E\), the \(E\)-span of $\left\{d_0\left(\xi(\varpi)^{4m}\prod_{j=0}^{m-1}\lambda_j\right)
\;\middle|\;
m\ge1
\right\}
$ is \(\widetilde E\). Hence \(\iota_0(\widetilde\tau_b)\subseteq W^\#\).
\end{proof}

Lemmas \ref{i_0 tilde} and \ref{propagation lemma} imply that $\iota_i(\widetilde{\tau}_b) \subseteq W^\#$ for all $i \in \mathbb{Z}$. 
Arguing as in the proof of Lemma~\ref{propagation lemma}, it follows 
that $\iota_i(\widetilde{D}_0^\#) \subseteq W^\#$ for all $i \in \mathbb{Z}$. 
Therefore $W^\# = \widetilde{D}_0^\#(\infty)$.
\end{proof}

\section{The irreducible \texorpdfstring{$G$}{G}-representation}

\begin{theorem}\label{main thm}
\begin{enumerate}
\item\label{1} There exists an $E[G]$-representation $\pi$ equipped with a commuting left $\widetilde{E}$-action, together with \(\widetilde E\)-linear injections
\(\widetilde D_0(\infty)\hookrightarrow\pi|_{KZ}\) and \(\widetilde D_1(\infty)\hookrightarrow\pi|_{N_G(I)}\)
of \(KZ\)- and \(N_G(I)\)-representations, respectively,
such that the diagram of \(IZ\)-representations
\[
\begin{tikzcd}
\widetilde{D}_0(\infty) \arrow[r, hook] 
& \pi|_{KZ} \\
\widetilde{D}_1(\infty) \arrow[r, hook] \arrow[u, hook]
& \pi|_{N_G(I)} \arrow[u, hook]
\end{tikzcd}
\]
commutes (the left vertical map is \(\phi\) from Lemma~\ref{D1R0}
and the right vertical map is the identity), and the map \(\widetilde D_0(\infty)\hookrightarrow\pi|_{KZ}\)
induces an isomorphism on \(E[K]\)-socles. Moreover, $\widetilde{D}_0(\infty)$ generates $\pi$ as an $E[G]$-representation.

\item Let $\pi$ be as in (1). For any irreducible smooth $E[K]$-representation $\sigma'$ and any $\varphi \in \mathrm{Hom}_{E[K]}(\sigma', \pi)$, we have
\[
\varphi \cdot \mathcal{H}(\sigma') 
= E \cdot \xi\left(
\begin{psmallmatrix}
\varpi&0\\
0&\varpi
\end{psmallmatrix}
\right)^\mathbb{Z} \cdot \varphi,
\]
where $\mathcal{H}(\sigma') := \mathrm{End}_{E[G]}(\mathrm{ind}_K^G \sigma')$ denotes the Hecke algebra.

In particular, if $\pi$ is irreducible, then the following are equivalent:
\begin{itemize}
\item $\pi$ admits a Hecke eigenvalue;
\item $\pi$ has a central character;
\item $\bigl(\begin{smallmatrix} \varpi & 0 \\ 0 & \varpi \end{smallmatrix}\bigr)$ acts on $\pi$ by a scalar in $E^\times$.
\end{itemize}

\item If $\lambda$ is as in Proposition~\ref{irred prop}, then any $\pi$ as in (1) is an irreducible smooth $E[G]$-representation.

\item Suppose that $\xi(\varpi)^4\lambda_0 \in E^\times$ and $\lambda_i \neq \lambda_0$ for all $i \neq 0$, and let $\pi$ be as in (1), irreducible as an $E[G]$-representation. Then $\mathrm{End}_{E[\mathrm{GL}_2(F)]}(\pi) \cong \widetilde{E}
$ as $E$-algebras.
\end{enumerate}
\end{theorem}

\begin{proof}

Fix an injective envelope $D_0^\#\hookrightarrow \operatorname{inj}_K(D_0^\#)
$ in the category of smooth \(E[K]\)-representations. Since
\(\operatorname{soc}_{E[K]} D_0^\#\) is essential in \(D_0^\#\),
\(\operatorname{inj}_K(D_0^\#)\) is also an injective envelope of
\(\operatorname{soc}_{E[K]} D_0^\#\). We have
\[
\operatorname{soc}_{E[K]} D_0^\#
=
\sigma\oplus\kappa^{[s]}
\oplus
\bigoplus_{\rho\in\{\sigma^{[s]},\tau,\tau^{[s]},\kappa\}}
(\rho_a\oplus\rho_b).
\]
Since this is a finite direct sum, the direct sum of the injective
envelopes of its summands is an injective envelope of
\(\operatorname{soc}_{E[K]} D_0^\#\). We therefore fix a decomposition
\[
\operatorname{inj}_K(D_0^\#)
\cong
\operatorname{inj}_K(\sigma)
\oplus
\operatorname{inj}_K(\kappa^{[s]})
\oplus
\bigoplus_{\rho\in\{\sigma^{[s]},\tau,\tau^{[s]},\kappa\}}
\left(J_{\rho,a}\oplus J_{\rho,b}\right),
\]
and identify the summands on the right with the corresponding
\(E[K]\)-subrepresentations of \(\operatorname{inj}_K(D_0^\#)\).
Here \(J_{\rho,a}\) and \(J_{\rho,b}\) are injective envelopes of \(\rho_a,\rho_b\subseteq D_0^\#\), respectively. That is, for each
\(\rho\in\{\sigma^{[s]},\tau,\tau^{[s]},\kappa\}\), we fix
\(K\)-isomorphisms
\begin{equation}\label{Jisos}
\operatorname{inj}_K(\rho)
\xrightarrow{\sim} J_{\rho,a}\quad \text{and}\quad \operatorname{inj}_K(\rho)
\xrightarrow{\sim} J_{\rho,b}
\end{equation}
extending the fixed isomorphisms \(\rho\xrightarrow{\sim}\rho_a\) and
\(\rho\xrightarrow{\sim}\rho_b\), respectively, in \eqref{socle-isomorphisms}. For each \(i\in\mathbb Z\), let \(
D_0^\#(i)\hookrightarrow
\operatorname{inj}_K(D_0^\#(i))
\) be the corresponding copy of this chosen injective envelope, together
with its corresponding decomposition. Set
\[
\Omega^\#
:=
\bigoplus_{i\in\mathbb Z}
\operatorname{inj}_K(D_0^\#(i))
\quad \text{and} \quad
\widetilde\Omega^\#
:=
\Omega^\#\otimes_E\widetilde E.
\]
The representation \(\widetilde\Omega^\#\) carries natural commuting
left and right \(\widetilde E\)-actions. The componentwise inclusions
induce a \(K\)-equivariant injection
\(
\widetilde D_0^\#(\infty)
\hookrightarrow
\widetilde\Omega^\#
\) commuting with both \(\widetilde E\)-actions. We extend the
\(K\)-action on \(\widetilde\Omega^\#\) to a \(KZ\)-action by letting
\(
\begin{psmallmatrix}
\varpi&0\\
0&\varpi
\end{psmallmatrix}^{\mathbb{Z}}
\) act through \(\xi\).

For \(H\in\{K,I\}\) and an \(E[H]\)-representation \(V\), write
\begin{equation}\label{injtilde}
    \widetilde{\operatorname{inj}}_H(V)
:=
\operatorname{inj}_H(V)\otimes_E\widetilde E.
\end{equation}
As before, for
\(\rho\in\{\sigma^{[s]},\tau,\tau^{[s]},\kappa\}\) and
\(v\in\widetilde{\operatorname{inj}}_K(\rho)\), write
\(v_a\) and \(v_b\) for its images in
\(\widetilde J_{\rho,a}:=J_{\rho,a}\otimes_E\widetilde E\)
and
\(\widetilde J_{\rho,b}:=J_{\rho,b}\otimes_E\widetilde E\),
respectively, under the scalar extensions of the fixed
isomorphisms in \eqref{Jisos}. Write \((v_a)_i\) and \((v_b)_i\)
for the corresponding elements in the \(i\)-th copy.
Analogously to \eqref{Del} and \eqref{DelLam}, we define \[
\Delta_{i,j}\bigl(\widetilde{\operatorname{inj}}_K(\rho)\bigr)
:=
\left\{
(v_a)_i-(v_b)_j
\;\middle|\;
v\in\widetilde{\operatorname{inj}}_K(\rho)
\right\}
\subseteq\widetilde\Omega^\#
\]
for \(\rho\in\{\sigma^{[s]},\tau^{[s]},\kappa\}\), and
\[
\Delta^{\lambda_i}_{i,i+1}
\bigl(\widetilde{\operatorname{inj}}_K(\tau)\bigr)
:=
\left\{
(v_a)_i-(v_b)_{i+1}\lambda_i
\;\middle|\;
v\in\widetilde{\operatorname{inj}}_K(\tau)
\right\}
\subseteq\widetilde\Omega^\#.
\]
Define
\[
S
:=
\bigoplus_{i\in\mathbb{Z}}
\Bigl(
\Delta_{i,i}
\bigl(\widetilde{\operatorname{inj}}_K(\sigma^{[s]})\bigr)
\oplus
\Delta^{\lambda_i}_{i,i+1}
\bigl(\widetilde{\operatorname{inj}}_K(\tau)\bigr)
\oplus
\Delta_{i,i-1}
\bigl(\widetilde{\operatorname{inj}}_K(\tau^{[s]})\bigr)
\oplus
\Delta_{i,i}
\bigl(\widetilde{\operatorname{inj}}_K(\kappa)\bigr)
\Bigr)
\subseteq
\widetilde{\Omega}^\#.
\]
Then \(S\) is an \(E[KZ]\)-subrepresentation of
\(\widetilde{\Omega}^\#\). Since each \(\lambda_i\) acts on the right,
\(S\) is also stable under the left
\(\widetilde E\)-action. We set
\[
\widetilde{\Omega}
:=
\widetilde{\Omega}^\#/S.
\]
Thus \(\widetilde{\Omega}\) is an \(E[KZ]\)-representation equipped
with a commuting left \(\widetilde E\)-action. Since the image of \(R\) under the componentwise
injection
\(
\widetilde D_0^\#(\infty)
\hookrightarrow
\widetilde{\Omega}^\#
\) is contained in \(S\), this injection induces a
\(KZ\)-equivariant, left \(\widetilde E\)-linear map $\widetilde D_0(\infty)
\rightarrow
\widetilde{\Omega}$.

\begin{lemma}\label{D0 injects Omega}
The intersection of the image of
\(\widetilde D_0^\#(\infty)\) with \(S\) inside
\(\widetilde{\Omega}^\#\) is precisely \(R\). In
particular, the induced map
\(
\widetilde D_0(\infty)
\rightarrow
\widetilde{\Omega}
\)
is injective.
\end{lemma}

\begin{proof}
Identify \(\widetilde D_0^\#(\infty)\) with its image in
\(\widetilde\Omega^\#\), and set $X:=\widetilde D_0^\#(\infty)\cap S$. Then \(\operatorname{soc}_{E[K]}X=R\). Indeed, as $R \subseteq \widetilde{D}_0^\#(\infty)$ and 
$R \subseteq S$, we have $R \subseteq X$, and as $R$ is semisimple, 
$R \subseteq \mathrm{soc}_{E[K]}(X)$. For the reverse inclusion, 
since $X \subseteq S$, we have 
$\mathrm{soc}_{E[K]}(X) \subseteq \mathrm{soc}_{E[K]}(S) = R$.

Let
\(q:\widetilde D_0^\#(\infty)\twoheadrightarrow
\widetilde D_0(\infty)\)
be the quotient map, and suppose that \(R\subsetneq X\).
Then \(q(X)\neq0\), so there exist a Serre weight \(\rho\) and an $E[K]$-linear injection $f:\rho\hookrightarrow
\operatorname{soc}_{E[K]}q(X)
\subseteq \operatorname{soc}_{E[K]}\widetilde D_0(\infty).
$ By Lemma~\ref{socle surjection}, \(f\) lifts to an $E[K]$-linear injection $f^\#:\rho\hookrightarrow
\operatorname{soc}_{E[K]}\widetilde D_0^\#(\infty)
$ such that \(q\circ f^\#=f\).
Since $f(\rho)\subseteq q(X)$ and \(\ker q=R\subseteq X\), we have
\(f^\#(\rho)\subseteq X\), and hence $f^\#(\rho)\subseteq\operatorname{soc}_{E[K]}X=R$. This implies \(q\circ f^\#=0\), contradicting \(f\neq0\).
Therefore \(X=R\).
\end{proof}

Since $\widetilde{D}_1(\infty)$ injects into $\widetilde{D}_0(\infty)|_{IZ}$ 
by Lemma~\ref{D1R0} and $\widetilde{D}_0(\infty)$ injects into 
$\widetilde{\Omega}$ by Lemma~\ref{D0 injects Omega}, we obtain a left \(\widetilde E\)-linear \(IZ\)-equivariant injection $\iota : \widetilde{D}_1(\infty) \hookrightarrow 
\widetilde{\Omega}|_{IZ}$. 
 
\begin{lemma}\label{N-action lemma}
There exists an action of \(N_G(I)\) 
extending the \(IZ\)-action on \(\widetilde{\Omega}\) such that the injection
\(
\iota:\widetilde D_1(\infty)
\hookrightarrow
\widetilde{\Omega}|_{IZ}
\)
is \(N_G(I)\)-equivariant.
\end{lemma}

\begin{proof}
Recall from the beginning of the proof of Theorem~\ref{main thm} that we
fixed an injective envelope \(D_0^\#\hookrightarrow\operatorname{inj}_K(D_0^\#)\) and, for \(\rho\in\{\sigma^{[s]},\tau,\tau^{[s]},\kappa\}\), \(E[K]\)-subrepresentations
\(J_{\rho,a}\) and \(J_{\rho,b}\), which are injective envelopes of
\(\rho_a\) and \(\rho_b\), respectively. Set
$$
\mathcal J:=
\operatorname{inj}_K(\sigma)
\oplus J_{\sigma^{[s]},a}
\oplus J_{\tau,b}
\oplus J_{\tau^{[s]},b}
\oplus J_{\kappa,a}
\oplus\operatorname{inj}_K(\kappa^{[s]})
\subseteq \operatorname{inj}_K(D_0^\#),
$$
and
$$
\mathcal J^{\mathrm c}:=
J_{\sigma^{[s]},b}
\oplus J_{\tau,a}
\oplus J_{\tau^{[s]},a}
\oplus J_{\kappa,b}.
$$
Thus $\operatorname{inj}_K(D_0^\#)=\mathcal J\oplus\mathcal J^{\mathrm c}$. For each \(i\in\mathbb Z\), let \(\mathcal J(i)\) and
\(\mathcal J^{\mathrm c}(i)\) be the corresponding subrepresentations of
\(\operatorname{inj}_K(D_0^\#(i))\), and set
$\mathcal J(\infty):=\bigoplus_{i\in\mathbb Z}\mathcal J(i)$. Recall that \(\widetilde V=V\otimes_E\widetilde E\) for any representation \(V\) over \(E\).

Consider the projection
$$
\widetilde\Omega^\#
\twoheadrightarrow
\bigoplus_{i\in\mathbb Z}\widetilde{\mathcal J^{\mathrm c}}(i).
$$
By the definition of \(S\), its restriction to
\(S\) is an isomorphism onto the target. Hence $\widetilde\Omega^\#
=\widetilde{\mathcal J}(\infty)\oplus S$. Therefore the quotient map
\(q:\widetilde\Omega^\#\twoheadrightarrow\widetilde\Omega\)
restricts to a left \(\widetilde E\)-linear \(KZ\)-isomorphism
$$
\varphi:\widetilde{\mathcal J}(\infty)
\xrightarrow{\sim}\widetilde\Omega.
$$
Set \(j:=\varphi^{-1}\circ\iota\), and let \(\mathscr C\) be the
category of smooth left \(\widetilde E[I]\)-modules.
By \cite[Corollary~6.4 and Lemma~6.8]{Bre07},
\begin{equation}\label{injective-restriction-to-I}
\operatorname{inj}_K(\rho)|_I
\cong
\operatorname{inj}_I\bigl(\chi(\rho)\bigr)
\oplus
\operatorname{inj}_I\bigl(\chi(\rho)^{(\Pi)}\bigr)
\end{equation}
for every generic Serre weight \(\rho\).
Thus each \(\mathcal J(i)|_I\) is injective in the category of smooth $E[I]$-modules. The proof of \cite[Lemma~5.29]{Pas13} applies to division algebras over \(E\), so extension of scalars to $\widetilde E$ preserves injectivity
for smooth \(E[I]\)-modules.
Moreover, since \(\mathscr C\) is locally Noetherian, direct
sums of injectives are injective by
\cite[Proposition~A.1.1(2)]{DEG}.
Hence \(\widetilde{\mathcal J}(\infty)|_I\) is an injective object in \(\mathscr C\).

Recall from the paragraph preceding \eqref{ABpsi} the
decomposition \(D_1=A\oplus B\), with a fixed isomorphism
\(\psi:A\xrightarrow{\sim}B^{(\Pi)}\).
Choose injective envelopes of \(A\) and \(B\) over \(E\),
and take their direct sum, giving an injective envelope
\(D_1\hookrightarrow\operatorname{inj}_I(D_1)
=\operatorname{inj}_I(A)\oplus\operatorname{inj}_I(B)\).
The isomorphism \(\psi\) extends to an \(E[I]\)-isomorphism
\(\operatorname{inj}_I(A)\xrightarrow{\sim}
\operatorname{inj}_I(B)^{(\Pi)}\).
Tensoring this map with \(\mathrm{id}_{\widetilde E}\) and
requiring \(\Pi^2\) to act by \(\xi(\varpi)\) defines an
\(N_G(I)\)-action on \(\widetilde{\operatorname{inj}}_I(D_1)\)
extending the given action on \(\widetilde D_1\). Here \(\widetilde{\operatorname{inj}}_I(D_1)
=\operatorname{inj}_I(D_1)\otimes_E\widetilde E\),
as defined in \eqref{injtilde}. For each \(i\in\mathbb Z\), let
\(D_1(i)\hookrightarrow\operatorname{inj}_I(D_1(i))\)
be a copy of this injective envelope, and give
\(\bigoplus_{i\in\mathbb Z}
\widetilde{\operatorname{inj}}_I(D_1(i))\)
the componentwise \(N_G(I)\)-action.
The same scalar-extension and direct-sum arguments as for \(\widetilde{\mathcal J}(\infty)|_I\) show that this representation is injective in \(\mathscr C\).
The inclusion
\(\widetilde D_1(\infty)\hookrightarrow
\bigoplus_{i\in\mathbb Z}
\widetilde{\operatorname{inj}}_I(D_1(i))\)
is essential, since it induces an isomorphism on
\(I_1\)-invariants.

By injectivity of \(\widetilde{\mathcal J}(\infty)|_I\),
\(j\) extends to a morphism
\(F:\bigoplus_{i\in\mathbb Z}
\widetilde{\operatorname{inj}}_I(D_1(i))
\to\widetilde{\mathcal J}(\infty)|_I\)
in \(\mathscr C\).
Since \(j\) is injective and \(\widetilde D_1(\infty)\) is
essential in the source, \(F\) is injective.
Its image \(\Omega_1:=\operatorname{im}(F)\) is therefore injective, so $\widetilde{\mathcal J}(\infty)|_I
=\Omega_1\oplus\Omega_2$ for some subrepresentation \(\Omega_2\).
As a direct summand of an injective object,
\(\Omega_2\) is also injective in \(\mathscr C\).

By construction,
\(D_1(i)^{I_1}
=(\operatorname{soc}_{E[K]}\mathcal J(i))^{I_1}\).
Since \(\varphi=q|_{\widetilde{\mathcal J}(\infty)}\),
\(j\) restricts to the coordinate inclusion on
\(\widetilde D_1(i)^{I_1}\).
Also,
\(\Omega_1^{I_1}=j(\widetilde D_1(\infty)^{I_1})\).
Consequently, we obtain
\[
\Omega_2^{I_1}
\cong
\frac{\widetilde{\mathcal J}(\infty)^{I_1}}{\Omega_1^{I_1}}
\cong
\bigoplus_{i\in\mathbb Z}
\frac{\widetilde{\mathcal J}(i)^{I_1}}
     {j\bigl(\widetilde D_1(i)^{I_1}\bigr)} \cong
\bigoplus_{i\in\mathbb Z}
\bigl(\widetilde D_1(i)^{I_1}\bigr)^{(\Pi)}
\cong
\bigl(\Omega_1^{I_1}\bigr)^{(\Pi)},
\]
where the third isomorphism follows from
\eqref{injective-restriction-to-I}.
Since \(\Omega_2\) and \(\Omega_1^{(\Pi)}\) are injective
envelopes of their socles in \(\mathscr C\), there is an
isomorphism \(\Omega_2\cong\Omega_1^{(\Pi)}\) in \(\mathscr C\).

Transport the \(N_G(I)\)-action to \(\Omega_1\) through \(F\),
and its \(\Pi\)-twist to \(\Omega_2\) through a fixed isomorphism \(\Omega_1^{(\Pi)}\xrightarrow{\sim}\Omega_2\)
in \(\mathscr C\).
The resulting action on
\(\widetilde{\mathcal J}(\infty)=\Omega_1\oplus\Omega_2\)
extends the given \(IZ\)-action and makes \(j\)
\(N_G(I)\)-equivariant.
Transporting this action through \(\varphi\) gives the required \(N_G(I)\)-action on \(\widetilde\Omega\), for which
\(\iota=\varphi\circ j\) is \(N_G(I)\)-equivariant.
\end{proof}

We now have the actions of $KZ$ and $N_G(I)$ on $\widetilde{\Omega}$ which agree on their intersection $IZ$. By \cite[Corollary~5.5.5]{Pas04}, these actions extend uniquely to an $E[G]$-action on $\widetilde{\Omega}$. Moreover, the natural left $\widetilde{E}$-action commutes with the $E[G]$-action, since it commutes with the actions of the subgroups $KZ$ and $N_G(I)$, which together generate $G$.

Let $\pi \subseteq \widetilde{\Omega}$ be the $E[G]$-subrepresentation generated by $\widetilde{D}_0(\infty)$. Then, by definition, we obtain the injections in part (\ref{1}) together with the commutative diagram, and $\widetilde{D}_0(\infty)$ generates $\pi$ as an $E[G]$-representation. Since $\widetilde{D}_0(\infty)$ is stable under the left $\widetilde{E}$-action and generates $\pi$ as an $E[G]$-representation, it follows that $\pi$ is also stable under $\widetilde{E}$. Thus, $\pi$ carries commuting left actions of $E[G]$ and $\widetilde{E}$, and the injection $\widetilde{D}_0(\infty) \hookrightarrow \pi$
is $\widetilde{E}$-linear and $KZ$-equivariant.

Consider the short exact sequence of left \(\widetilde E[K]\)-modules
\[
0\longrightarrow S
\longrightarrow \widetilde\Omega^\#
\xrightarrow{q}
\widetilde\Omega
\longrightarrow 0.
\]
By the proof of Lemma~\ref{N-action lemma}, the quotient map
\(q:\widetilde\Omega^\#\twoheadrightarrow\widetilde\Omega\)
restricts to a left \(\widetilde E\)-linear \(KZ\)-isomorphism
\(
\varphi:
\widetilde{\mathcal J}(\infty)
\xrightarrow{\sim}
\widetilde\Omega.
\) Hence
\(
\widetilde\Omega
\xrightarrow{\varphi^{-1}}
\widetilde{\mathcal J}(\infty)
\hookrightarrow
\widetilde\Omega^\#
\) is a section of \(q\). Therefore, the above short exact sequence splits. Consequently, $\mathrm{soc}_{E[K]}\widetilde{\Omega}^\# \rightarrow \mathrm{soc}_{E[K]}\widetilde{\Omega}$ is surjective.

Furthermore, since the inclusion $\widetilde{D}_0^\#(\infty) \hookrightarrow \widetilde{\Omega}^\#$ is a direct sum of essential extensions, it induces an isomorphism on $E[K]$-socles. Combining these facts with the left-exactness of the socle functor applied to the injection $\widetilde{D}_0(\infty) \hookrightarrow \widetilde{\Omega}$, we obtain the following commutative diagram:
\[
\begin{tikzcd}[row sep=large, column sep=large]
\mathrm{soc}_{E[K]}\widetilde{\Omega}^\# \arrow[r, two heads] & \mathrm{soc}_{E[K]}\widetilde{\Omega} \\
\mathrm{soc}_{E[K]}\widetilde{D}_0^\#(\infty) \arrow[u, "\sim" sloped] \arrow[r] & \mathrm{soc}_{E[K]}\widetilde{D}_0(\infty) \arrow[u, hook]
\end{tikzcd}.
\]
This forces the inclusion $\mathrm{soc}_{E[K]}\widetilde{D}_0(\infty) \hookrightarrow \mathrm{soc}_{E[K]}\widetilde{\Omega}$ to be an isomorphism. Finally, since $\pi \subseteq \widetilde{\Omega}$ is the $E[G]$-subrepresentation generated by $\widetilde{D}_0(\infty)$, the chain of inclusions $\widetilde{D}_0(\infty) \subseteq \pi \subseteq \widetilde{\Omega}$ ensures that $\widetilde{D}_0(\infty) \hookrightarrow \pi$ induces the desired isomorphism on $E[K]$-socles. This completes the proof of part~(\ref{1}).

Next, we prove \textup{(2)}. We may assume that \(\varphi\neq0\).
By part~\textup{(1)}, \(\sigma'\) is isomorphic to one of
\(\sigma\), \(\sigma^{[s]}\), \(\tau\), \(\tau^{[s]}\), \(\kappa\), and \(\kappa^{[s]}\).
The Hecke algebra \(\mathcal H(\sigma')\) is generated over \(E\)
by the Hecke operator \(T\), unique up to scalar, supported on
\(K\bigl(\begin{smallmatrix}1&0\\0&\varpi\end{smallmatrix}\bigr)K\),
and the central operators acting by
\(\xi(\varpi)^{\pm1}\) \cite{Her11b}.
We show that \(\varphi\cdot T=0\).

Suppose first that \(\sigma'=\tau\), and let
\(\Phi:\operatorname{ind}_K^G\tau\to\pi\)
be the \(E[G]\)-map corresponding to \(\varphi\) under
Frobenius reciprocity.
Choose \(0\neq v_0\in\tau^{I_1}\) and put \(v:=\varphi(v_0)\).
As in the proof of \cite[Lemma~4.1]{GLS23}, we have
\(
(\Phi\circ T)(\tau)\subseteq \xi(\varpi)^{-1}(K\cdot\Pi v) \cong\tau^{[s]}\).
Since \(\operatorname{Hom}_{E[K]}(\tau,\tau^{[s]})=0\),
we obtain \((\Phi\circ T)|_\tau=0\), and hence
\(\varphi\cdot T=0\).
The other cases are identical.
This proves the first assertion of~\textup{(2)}.
The equivalence of the three conditions follows as in the proof
of \cite[Theorem~3.1(2)]{Le24}.

We now prove~(3). Let $\pi' \subseteq \pi$ be a nonzero $E[G]$-subrepresentation, and set $W := \pi'|_{KZ} \cap \widetilde{D}_0(\infty).
$ Then $W$ contains the $E[K]$-socle of $\pi'|_K$ and is therefore nonzero.

Let $W^\#$ be the preimage of $W$ under the quotient map 
$\widetilde{D}_0^\#(\infty) \twoheadrightarrow \widetilde{D}_0(\infty)$. 
Then $W^\#$ is an $E[KZ]$-subrepresentation of $\widetilde{D}_0^\#(\infty)$ such that $R\subsetneq W^\#$. Moreover, since $\pi'$ is stable under the action of $G$, it is in particular stable under $\Pi$, and it follows that $W^\# \cap \widetilde{D}_1(\infty)
$ is stable under the action of $\Pi$. We may therefore apply Proposition \ref{irred prop} to $W^\#$, which shows that $W^\# = \widetilde{D}_0^\#(\infty).$ Passing to the quotient, we deduce that $W = \widetilde{D}_0(\infty)$.
It follows that $\pi'|_{KZ}$ contains $\widetilde{D}_0(\infty)$, and hence $\pi'$ contains the $E[G]$-subrepresentation generated by $\widetilde{D}_0(\infty)$, which is $\pi$. Thus, $\pi' = \pi$, and $\pi$ is irreducible.

Finally, we prove (4). The left \(\widetilde E\)-action on \(\pi\) gives an
injective homomorphism \(\widetilde E\hookrightarrow \End_{E[G]}(\pi)\) of
\(E\)-algebras. We show that it is an isomorphism, following
\cite[Theorem~3.1(4)]{Le24}. Let
\[
\jmath:\tau\xrightarrow{\sim}\tau_b
\hookrightarrow \widetilde{\tau}_b
\xrightarrow{\ \iota_0\ } \M
\twoheadrightarrow \widetilde D_0(\infty)
\hookrightarrow \pi
\]
be the natural map, where the first map is the fixed isomorphism in
\eqref{socle-isomorphisms} and \(\iota_0\) is restricted
to \(\widetilde\tau_b\subseteq\widetilde D_0^\#\). Since \(\pi\) is irreducible, precomposition with
\(\jmath\) gives an injection \(
\End_{E[G]}(\pi)\hookrightarrow \Hom_{E[K]}(\tau,\pi)
\). It is enough to show that its image is contained in
\(\Hom_{E[K]}(\tau,(\iota_0(\widetilde{\tau}_b)+R)/R)\cong
\widetilde E\). Granting this, the composition
\(
\widetilde E\hookrightarrow \End_{E[G]}(\pi)\hookrightarrow \widetilde E
\) is the identity by the \(\widetilde E\)-linearity in (1), and we are done.

Write
\[
R_{\sigma^{[s]}}
:=\bigoplus_{i\in\mathbb Z}
\Delta_{i,i}(\widetilde{\sigma}^{[s]}),\,
R_{\tau}:=\bigoplus_{i\in\mathbb Z}
\Delta_{i,i+1}^{\lambda_i}(\widetilde{\tau}),\,
R_{\tau^{[s]}}:=\bigoplus_{i\in\mathbb Z}\Delta_{i,i-1}(\widetilde{\tau}^{[s]}),\,\text{and}\,
R_{\kappa}:=\bigoplus_{i\in\mathbb Z}
\Delta_{i,i}(\widetilde{\kappa}).
\]
Thus $R=
R_{\sigma^{[s]}}
\oplus
R_\tau
\oplus
R_{\tau^{[s]}}
\oplus
R_\kappa$. The subrepresentations
\(R_{\sigma^{[s]}},R_{\tau^{[s]}}\), and
\(R_\kappa\) are invariant under translations of the index set
\(\mathbb Z\). The term \(R_\tau\) is not, in general, because it
contains the scalars \(\lambda_i\).

Choose a nonzero vector \(v_b\in\tau_b^{I_1}\). Applying \(\Pi\) to
\((v_b)_i\) gives \((\Pi v_b)_i\in
\iota_i(\widetilde{\tau}_b^{[s]})^{I_1}\), and \(K\cdot(\Pi v_b)_i= \iota_i(\tau_b^{[s]})\); in particular, this \(K\)-span contains the $I$-submodule \(\iota_i(V_2(\tau_b^{[s]}))\).
The proof of the first containment in
Lemma~\ref{propagation lemma}, using only the
\(R_{\sigma^{[s]}}\)-relations, gives an element
\(S_1\in E[G]\) such that
\begin{equation}\label{S1action}
S_1\cdot ((v_b)_i+R)
=
\xi(\varpi)^2 (v_a)_i+R
\qquad\text{for all } i\in\mathbb Z.
\end{equation}
The same \(S_1\) works for every \(i\), because the \(K\)- and \(\Pi\)-actions on $\widetilde D_0^\#(\infty)$ and $\widetilde D_1(\infty)$, respectively, commute with shifting indices (see \eqref{Pi action}), and
\(R_{\sigma^{[s]}}\) is invariant under shifting. Now, use the \(\tau\)-gluing relation. Since
\(R_\tau\) contains
\(\Delta^{\lambda_i}_{i,i+1}(\widetilde\tau)\), we have
\begin{align}
(v_a)_i+R
&=(v_b)_{i+1}\lambda_i+R \notag\\
&=\lambda_i(v_b)_{i+1}+R
\qquad\text{for all }i\in\mathbb Z.
\label{tau-gluing}
\end{align}
Here, the second equality uses that
\(v_b\in\tau_b\) is identified with
\(v_b\otimes 1\in\widetilde\tau_b\), so left and right
multiplication by elements of \(\widetilde E\) agree on \(v_b\). Combining \eqref{S1action} and \eqref{tau-gluing}, we get
\begin{equation}\label{finalS1action}
S_1\cdot ((v_b)_i+R)
=
\xi(\varpi)^2 \lambda_i (v_b)_{i+1} +R
\qquad\text{for all } i\in\mathbb Z.
\end{equation}

Similarly, the proof of the second containment in
Lemma~\ref{propagation lemma} uses only the
\(R_\kappa\)- and
\(R_{\tau^{[s]}}\)-relations, which are invariant under
shifting indices, together with the \(K\)- and \(\Pi\)-actions, which commute
with shifting indices. Thus, there is an element \(S_2\in E[G]\) such
that
\begin{equation}
S_2\cdot\bigl((v_b)_{i+1}+R\bigr)
=
\xi(\varpi)^2(v_b)_i+R
\qquad\text{for all }i\in\mathbb Z.
\label{S2action}
\end{equation}

Set \(S_3:=S_2S_1\). Combining \eqref{finalS1action} and \eqref{S2action}, we obtain
\begin{equation}\label{S-action}
S_3\cdot\bigl((v_b)_i+R\bigr)
=
\xi(\varpi)^4\lambda_i(v_b)_i+R
\qquad\text{for all }i\in\mathbb Z.
\end{equation}

Now, let \(\phi\in\End_{E[G]}(\pi)\). Since \(\phi\) commutes with the
\(K\)-action, the vector \(\phi((v_b)_0+R)\) is \(I_1\)-invariant
and lies in the \(\tau\)-isotypic part of \(\soc_{E[K]}(\pi)\). In
\(\widetilde D_0(\infty)=\widetilde D_0^\#(\infty)/R\), the
relations in \(R_\tau\) identify
\(\iota_i(\widetilde\tau_a)+R\) with
\(\iota_{i+1}(\widetilde\tau_b)+R\) as subspaces. Hence, the subspace of $I_1$-invariants in the \(\tau\)-isotypic part of \(\soc_{E[K]}(\widetilde D_0(\infty))\) is represented by
\(
\bigoplus_{i\in\mathbb Z}\iota_i(\widetilde\tau_b^{I_1})+R.
\)
Thus,
\[
\phi\bigl((v_b)_0+R\bigr)
=
\sum_i c_i (v_b)_i+R
\]
for some \((c_i) \in
\bigoplus_{i \in \mathbb Z} \widetilde E\). Since \(\phi\) commutes with \(S_3\), equation~\eqref{S-action} gives
\begin{align*}
\sum_i c_i\,\xi(\varpi)^4\lambda_i\,(v_b)_i+R
&=
S_3\bigl(\phi((v_b)_0+R)\bigr) \\
&=
\phi\bigl(S_3((v_b)_0+R)\bigr) \\
&=
\phi\bigl(\xi(\varpi)^4\lambda_0 (v_b)_0+R\bigr) \\
&=
\sum_i c_i\,\xi(\varpi)^4\lambda_0\,(v_b)_i+R.
\end{align*}
Here we use the fact that the left \(\widetilde E\)-action commutes with \(E[G]\), and
that \(\xi(\varpi)^4\lambda_0\in E^\times\subseteq\widetilde C^\times\)
commutes with every \(c_i\in\widetilde E\). Thus,
\[
\sum_i c_i\,\xi(\varpi)^4(\lambda_i-\lambda_0)(v_b)_i\in R.
\]

We claim that
\(
R\cap
\bigoplus_{i\in\mathbb Z}\iota_i(\widetilde\tau_b^{I_1})
=0.
\) Indeed, the projection
\[
\operatorname{soc}_{E[K]}\widetilde D_0^\#(\infty)
\twoheadrightarrow
\bigoplus_{i\in\mathbb Z}
\left(
\iota_i(\widetilde\sigma_a^{[s]})
\oplus
\iota_i(\widetilde\tau_a)
\oplus
\iota_i(\widetilde\tau_a^{[s]})
\oplus
\iota_i(\widetilde\kappa_a)
\right)
\]
is injective on \(R\), by the definition of
\(R\), whereas it vanishes on
\(
\bigoplus_i\iota_i(\widetilde\tau_b^{I_1})
\). Thus $\sum_i c_i\,\xi(\varpi)^4(\lambda_i-\lambda_0)(v_b)_i=0$, and hence $c_i\,\xi(\varpi)^4(\lambda_i-\lambda_0)=0$ for all $i \in \mathbb Z$. Since \(\xi(\varpi)^4\in\widetilde C^\times\), \(\widetilde E\) is a division
algebra, and \(\lambda_i\neq\lambda_0\) for \(i\neq0\), we have \(c_i=0\) for
every \(i\neq0\). Hence,
\[
\phi\bigl((v_b)_0+R\bigr)
=
c_0(v_b)_0+R
\in
(\iota_0(\widetilde\tau_b^{I_1})+R)/R.
\]
Since \(v_b\) generates \(\tau_b\) as an \(E[K]\)-module and \(\phi\) commutes
with \(K\), this gives
\(
\phi\circ\jmath(\tau)
\subseteq
(\iota_0(\widetilde\tau_b)+R)/R.
\) Thus the restriction map factors through
\(\Hom_{E[K]}(\tau,(\iota_0(\widetilde\tau_b)+R)/R)\cong
\widetilde{E}\), as required.
\end{proof}

\section{The passage to \texorpdfstring{\(\mathrm{GL}_n\)}{GLn}}

We now pass from \(\mathrm{GL}_2\) to \(\mathrm{GL}_n\). We use
\cite[Proposition~3.3]{Le24}, which applies to arbitrary nonarchimedean local
fields and standard lower triangular parabolic subgroups.

\begin{theorem}\label{GLn theorem}
Let \(p>7\) be a prime and \(n>1\) be an integer, or let \(n=2\) and \(p=7\). Let \(F\neq\mathbb Q_p\) be a
nonarchimedean local field with residue field \(\mathbb F_p\). Let \(E\) be
an extension of \(\mathbb F_p\), and let \(\widetilde E\) be a division
algebra over \(E\) of countable \(E\)-dimension, with center
\(\widetilde C\). Then there exists a nonadmissible irreducible smooth
\(E[\mathrm{GL}_n(F)]\)-representation \(\rho\) such that \(\End_{E[\mathrm{GL}_n(F)]}(\rho)\cong \widetilde E
\) as \(E\)-algebras. Moreover, one may choose \(\rho\)
\begin{enumerate}
\item to admit a Hecke eigenvalue and hence to have a central character;
\item to have no central character, if \(E\subsetneq\widetilde C\);
\item to have a central character but no Hecke eigenvalue, if
\(n>3\) and \(E\subsetneq\widetilde C\).
\end{enumerate}
\end{theorem}

\begin{proof}

Given a \(\widetilde C\)-valued character \(\xi\) of
\(
\begin{psmallmatrix}
\varpi&0\\
0&\varpi
\end{psmallmatrix}^{\mathbb{Z}}
\), choose
\(\lambda=(\lambda_i)_{i\in\mathbb Z}\in
\prod_{i\in\mathbb Z}\widetilde E^\times\)
satisfying the hypotheses of Proposition~\ref{irred prop} as follows. Let
\((\mu_i)_{i\geq0}\) be an \(E\)-spanning sequence in
\(\widetilde E^\times\) such that \(\mu_0=1\) and
\(\mu_i\neq\mu_{i-1}\) for \(i>0\), and set
\[
\lambda_0:=\xi(\varpi)^{-4},\quad
\lambda_i:=\xi(\varpi)^{-4}\mu_{i-1}^{-1}\mu_i\quad\text{for all}\quad i>0,
\quad\text{and}\quad
\lambda_i:=-\lambda_0\quad\text{for all}\quad i<0.
\]
Then \(\xi(\varpi)^4\lambda_0=1\in E^\times\),
\(\lambda_i\neq\lambda_0\) for \(i\neq0\), and $\prod_{j=0}^{m-1}\xi(\varpi)^4\lambda_j=\mu_{m-1}$ for all $m\geq1$. Thus, the \(E\)-span of
\(\{\prod_{j=0}^{m-1}\xi(\varpi)^4\lambda_j\mid m\geq1\}\)
is \(\widetilde E\). Let \(\pi\) be the representation supplied by
Theorem~\ref{main thm} for these choices of \(\xi\) and \(\lambda\). Then
\(\pi\) is irreducible by Theorem~\ref{main thm}\textup{(3)}, and
\(
\End_{E[\GL_2(F)]}(\pi)\cong\widetilde E
\)
by Theorem~\ref{main thm}\textup{(4)}. Moreover, \(\pi\) is nonadmissible, since
\(\operatorname{soc}_{E[K]}\pi\cong
\operatorname{soc}_{E[K]}\widetilde D_0(\infty)\) contains infinitely
many copies of a fixed Serre weight.

If \(n=2\), for \textup{(1)} choose \(\xi\) with
\(\xi(\varpi)\in E^\times\), and for \textup{(2)}, when
\(E\subsetneq\widetilde C\), choose
\(\xi(\varpi)\in\widetilde C^\times\setminus E^\times\).
The assertions then follow from Theorem~\ref{main thm}\textup{(2)}. We therefore assume
that \(n>2\) and $p>7$.

Let \(P\) be the lower block triangular parabolic subgroup of
\(\GL_n(F)\) with block sizes \(2,1,\ldots,1\), with Levi quotient
\(
M\cong\GL_2(F)\times(F^\times)^{n-2}.
\)
Let \(\omega:\mathcal O_F^\times\to E^\times\) be the character
\(u\mapsto\overline u\). Choose an integer \(c\leq-1\) such that
\(c\not\equiv y\pmod{p-1}\) for every weight \(F(x,y)\) among
\(\sigma,\sigma^{[s]},\tau,\tau^{[s]},\kappa,\kappa^{[s]}\).
Such a choice is possible since \(p-1>6\). For a choice of \(\xi\) and the corresponding representation
\(\pi\) constructed above, choose a smooth character
\(\chi_0:(F^\times)^{n-2}\to\widetilde C^\times\) whose restriction
to \((\mathcal O_F^\times)^{n-2}\) is
\[
\bigotimes_{j=0}^{n-3}\omega^{c-j};
\]
its values on uniformizers will be specified below. Set
\[
\rho:=\Ind_P^{\GL_n(F)}
(\pi\otimes_{\widetilde C}\chi_0).
\]

Every \(\GL_n(\mathcal O_F)\)-weight occurring in \(\rho\)
is \(M\)-regular. Indeed, let \(F(x_1,\ldots,x_n)\) be such
a weight, with \(0\leq x_i-x_{i+1}\leq p-1\).
By Frobenius reciprocity and the description of
\(\operatorname{soc}_{E[K]}\pi\), the weight \(F(x_1,x_2)\)
is isomorphic to one of \(\sigma\), \(\sigma^{[s]}\), \(\tau\), \(\tau^{[s]}\), \(\kappa\), and \(\kappa^{[s]}\). We also get $x_j\equiv c+3-j\pmod{p-1}
$ for \(3\leq j\leq n\). Our choice of \(c\) gives \(x_2\not\equiv x_3\pmod{p-1}\),
while
\(x_j-x_{j+1}\equiv1\pmod{p-1}\) for \(3\leq j<n\).
Thus, \(x_i-x_{i+1}>0\) for every \(1\leq i<n\), so the
weight is \(M\)-regular by
\cite[Definition~2.4]{Her11a}.

For \textup{(1)}, choose \(\xi\) with
\(\xi(\varpi)\in E^\times\), construct \(\lambda\) and \(\pi\) as
above, and choose \(\chi_0\) to be \(E^\times\)-valued. The
Hecke-compatibility hypothesis of \cite[Proposition~3.3]{Le24} is
verified exactly as in the proof of \cite[Theorem~3.5]{Le24}, using
Theorem~\ref{main thm}\textup{(2)} in place of \cite[Theorem~3.1(2)]{Le24}. Hence
\(\rho\) is irreducible by \cite[Proposition~3.3]{Le24}. It is
nonadmissible because \(\operatorname{soc}_{E[K]}\pi\) contains
infinitely many copies of a fixed Serre weight, and the evaluation
isomorphism in \cite[Proposition~3.3]{Le24} identifies these with
infinitely many copies of the corresponding \(M\)-regular
\(\GL_n(\mathcal O_F)\)-weight in \(\rho\). Since
\(\xi(\varpi)\in E^\times\), Theorem~\ref{main thm}\textup{(2)} gives a Hecke eigenvalue
for \(\pi\), and hence \(\rho\) admits a Hecke eigenvalue by
\cite[Proposition~3.3]{Le24}. In particular, \(\rho\) has a central
character.

For \textup{(2)}, assume \(E\subsetneq\widetilde C\), choose
\(\xi\) with
\(\xi(\varpi)\in\widetilde C^\times\setminus E^\times\), construct
\(\lambda\) and \(\pi\) as above, and again choose \(\chi_0\) to be
\(E^\times\)-valued. The preceding argument shows that \(\rho\) is
irreducible and nonadmissible. The element \(\varpi I_n\) acts on
\(\rho\) by $\xi(\varpi)\chi_0(\varpi,\ldots,\varpi)
\notin E^\times,$ so \(\rho\) has no central character.

For \textup{(3)}, assume \(n>3\) and
\(E\subsetneq\widetilde C\). Choose
\(X\in\widetilde C^\times\setminus E^\times\), choose \(\xi\) with
\(\xi(\varpi)=X\), and construct \(\lambda\) and \(\pi\) as above.
Write
\(\chi_0=\chi_1\otimes\cdots\otimes\chi_{n-2}\), and, without changing the prescribed restrictions to units, choose
\[
\chi_1(\varpi)=X^{-2},\quad
\chi_2(\varpi)=X,\quad\text{and}\quad
\chi_j(\varpi)\in E^\times\quad \text{for all}\quad 2<j\leq n-2.
\]
As in the last paragraph of the proof of
\cite[Theorem~3.5]{Le24}, using Theorem~\ref{main thm}\textup{(2)} in place of
\cite[Theorem~3.1(2)]{Le24}, the Hecke-compatibility hypothesis of
\cite[Proposition~3.3]{Le24} holds. Thus, \(\rho\) is irreducible. It is nonadmissible by the same
weight-multiplicity argument as above. Moreover, \(\varpi I_n\)
acts on \(\rho\) by a scalar in \(E^\times\), whereas
\(\pi\otimes_{\widetilde C}\chi_0\) has no Hecke eigenvalue.
Therefore \(\rho\) has a central character but no Hecke eigenvalue,
by \cite[Proposition~3.3]{Le24}.

Finally, in each case full faithfulness of parabolic induction gives
\[
\End_{E[\GL_n(F)]}(\rho)
\cong
\End_{E[M]}(\pi\otimes_{\widetilde C}\chi_0)
\cong
\End_{E[\GL_2(F)]}(\pi)
\cong
\widetilde E,
\]
where the first isomorphism follows from \cite[Theorem~5.3(1)]{Vig16} and the last from Theorem~\ref{main thm}\textup{(4)}.
This completes the proof.
\end{proof}

\begin{proof}[Proof of Theorem~\ref{thm:main}]
If the residue field of \(F\) is different from \(\mathbb F_p\), the
result follows from \cite[Theorem~1.1]{Le24}. It remains to consider
the case where the residue field of \(F\) is \(\mathbb F_p\). The second assertion follows directly from
Theorem~\ref{GLn theorem}. For the first assertion, take
\(\widetilde E=\widetilde C\) to be a proper field extension of \(E\)
of countable degree, for instance \(E(t)\), and apply
Theorem~\ref{GLn theorem}.
\end{proof}

\end{document}